\documentclass[a4paper,twoside,11pt]{article}
\usepackage{amsfonts,amssymb,mathrsfs,amsmath}
\usepackage[dvips]{color}
\usepackage{graphicx}
\title{Towards Resolving Keller's Cube Tiling Conjecture in Dimension Seven}
\date{}
\author{Andrzej P. Kisielewicz\\
\\
{\small Wydzia{\l} Matematyki, Informatyki i Ekonometrii, Uniwersytet Zielonog\'orski}\\
{\small ul. Z. Szafrana 4a, 65-516 Zielona G\'ora, Poland}\\
{\small A.Kisielewicz@wmie.uz.zgora.pl}\\
}

\numberwithin{equation}{section}
\newtheorem{pr}{\sc Proposition}
\newtheorem{ex}[pr]{\sc Example}
\newtheorem{lemat}[pr]{\sc Lemma}

\newtheorem{tw}[pr]{\sc Theorem}
\newtheorem{wn}[pr]{\sc Corollary}
\newtheorem{df}{\sc Definition}

\newtheorem{uw}{\sc Remark}
\newtheorem{uwi}[uw]{\sc Remarks}

\newtheorem{nap}{\sc Example }

\newtheorem{nps}[nap]{\sc Examples}

\def\ka #1{\mathscr{#1}}
\def\kal #1 #2{\mathscr{#1}^{#2}}
\def\proof{\noindent \textit{Proof.\,\,\,}}

\def\zet{\mathbb{Z}}

\def\er{\mathbb{R}}

\def\te{\mathbb{T}}
\def\en{\mathbb{N}}

\def\iver #1{\mbox{\tt [} #1 \mbox{\tt]}}

\begin{document}

\numberwithin{pr}{section}
\numberwithin{uw}{section}

\maketitle
\begin{abstract}
A cube tiling of $\er^d$ is a  family of pairwise disjoint cubes $[0,1)^d+T=\{[0,1)^d+t\colon t\in T\}$ such that $\bigcup_{t\in T}([0,1)^d+t)=\er^d$. Two cubes $[0,1)^d+t$, $[0,1)^d+s$ are called a twin pair if $|t_j-s_j|=1$ for some $j\in [d]=\{1,\ldots, d\}$ and $t_i=s_i$ for every $i\in [d]\setminus \{j\}$. 
In $1930$, Keller conjectured that in every cube tiling of $\er^d$ there is a twin pair. Keller's conjecture is true for dimensions $d\leq 6$ and false for all dimensions $d\geq 8$. For $d=7$ the conjecture is still open. Let $x\in \er^d$, $i\in [d]$, and let $L(T,x,i)$ be the set of all $i$th coordinates $t_i$ of vectors $t\in T$ such that $([0,1)^d+t)\cap ([0,1]^d+x)\neq \emptyset$ and $t_i\leq x_i$. Let $r^-(T)=\min_{x\in \er^d}\; \max_{1\leq i\leq d}|L(T,x,i)|$ and $r^+(T)=\max_{x\in \er^d}\; \max_{1\leq i\leq d}|L(T,x,i)|$.
It is known that if $r^-(T)\leq 2$ or  $r^+(T)\geq 5$, then Keller's conjecture is true for $d=7$. In the paper we show that it is also true for $d=7$ if $r^+(T)=4$. Thus, if $[0,1)^7+T$ is a counterexample to Keller's conjecture, then $r^+(T)=3$, which is the last unsolved case of Keller's conjecture. Additionally, a new proof of Keller's conjecture in dimensions $d\leq 6$ is given.



\end{abstract}
\section{Introduction}
A {\it cube tiling} of $\er^d$ is a  family of pairwise disjoint cubes $[0,1)^d+T=\{[0,1)^d+t\colon t\in T\}$ such that $\bigcup_{t\in T}([0,1)^d+t)=\er^d$. Two cubes $[0,1)^d+t$, $[0,1)^d+s$ are called a {\it twin pair} if $|t_j-s_j|=1$
for some $j\in [d]=\{1,\ldots, d\}$ and $t_i=s_i$ for every $i\in [d]\setminus \{j\}$. In $1930$, Keller \cite{Ke1} conjectured that in every cube tiling of $\er^d$ there is a twin pair. In $1940$, Perron \cite{P} proved that Keller's conjecture is true for all dimensions $d\leq 6$ (see also \cite{LP1}).  
In 1992, Lagarias and Shor \cite{LS1}, using ideas from  Corr\'adi's and Szab\'o's papers \cite{CS2,Sz2}, constructed a cube tiling of $\er^{10}$ which does not contain a twin pair  and thereby refuted  Keller's cube tiling conjecture. In $2002$, Mackey \cite{M} gave a counterexample to Keller's conjecture in dimension eight, which also shows that this conjecture is false in dimension nine. For $d=7$ Keller's conjecture is still open.  

Let  $[0,1)^d+T$ be a cube tiling, $x\in \er^d$ and $i\in [d]$, and let $L(T,x,i)$ be the set of all $i$th coordinates $t_i$ of vectors $t\in T$ such that $([0,1)^d+t)\cap ([0,1]^d+x)\neq \emptyset$ and $t_i\leq x_i$ (Figure 1). 
For every $x\in \er^d$ and $i\in [d]$ the set $L(T,x,i)$ contains at most $2^{d-1}$ elements.

\smallskip
{\center
\includegraphics[width=6cm]{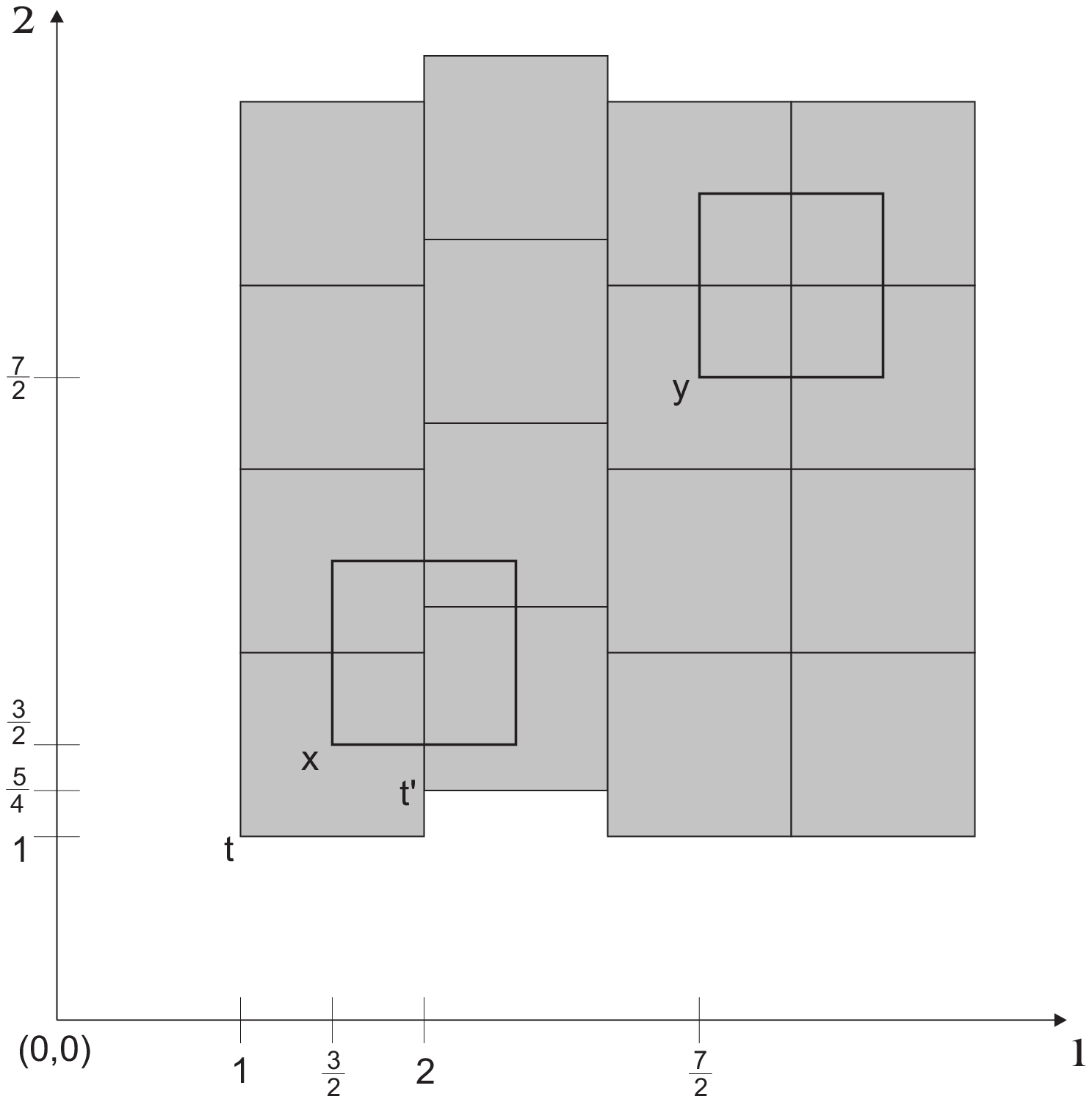}\\
}

\noindent
{\footnotesize 
Figure 1: A portion of a cube tiling $[0,1)^2+T$ of $\er^2$. The number of elements in $L(T,x,i)$ depends on the position of $x\in \er^2$. In the picture we have $x=(\frac{3}{2},\frac{3}{2})$, $y=(\frac{7}{2},\frac{7}{2})$, $t=(1,1)$ and $t'=(2,\frac{5}{4})$. Then $L(T,x,1)=\{1\}$, $L(T,x,2)=\{1,\frac{5}{4}\}$ and $L(T,y,1)=L(T,y,2)=\{3\}$. This portion of the tiling  $[0,1)^2+T$ shows that $r^-(T)=1$ and $r^+(T)=2$.
  }

\smallskip
Let 
\begin{equation}
\label{re}
r^-(T)=\min_{x\in \er^d}\; \max_{1\leq i\leq d}|L(T,x,i)|\;\;\;\; {\rm and}\;\;\;\; r^+(T)=\max_{x\in \er^d}\; \max_{1\leq i\leq d}|L(T,x,i)|.
\end{equation}

In 2010, Debroni et al. \cite{De} computed that Keller's conjecture is true for all cube tilings $[0,1)^7+T$ of $\er^7$ such that  $T\subset \frac{1}{2}\zet^7$. This result shows that Keller's conjecture is true for cube tilings of $\er^7$ with $r^-(T)\leq 2$ (\cite{Kis}). 
In \cite{Kis,Kis22} we showed that Keller's conjecture is true for cube tilings $[0,1)^7+T$ of $\er^7$ for which $r^+(T)\geq 5$.
In the presented paper we prove that 

\begin{tw}
\label{keli1}
Keller's conjecture is true for all cube tilings $[0,1)^7+T$ of $\er^7$ for which $r^+(T)=4$.
\end{tw}

Thus, the above theorem resolves the penultimate case of Keller's conjecture in dimension seven. To complete resolution of this conjecture remains to resolve the last case $r^+(T)=3$. (Corollary \ref{kon3}).

Our proof of Theorem \ref{keli1} is based on a structural result dealing with two systems of abstract words  having $n\in \{13,...,16\}$ words each (Theorem \ref{ke2}). It can be interpreted by means of systems of cubes in the flat torus  $\te^d=\{(x_1,\ldots ,x_d)({\rm mod} 2):(x_1,\ldots ,x_d)\in \er^d\}$ as follows:
A set $F\subset \te^d$ is called a {\it polycube} if $F$ has a { \it tiling} by translates of the unit cube, that is, there is a family of pairwise disjoint translates of the unit cube $\ka F=[0,1)^d+T$, $T\subset \te^d$, such that $\bigcup_{t\in T}[0,1)^d+t=F$. 
The case $r^+(T)=4$ is reduced to the following task: {\it For $d=6$ determine all polycubes $F\subset \te^d$ which have at least two twin pair free cube tilings $\ka F$ and $\ka G$ such that $\ka F\cap \ka G=\emptyset$ and $|\ka F|=|\ka G|=n$, where $n\in \{13,...,16\}$}.  As we shall show (Theorem \ref{ke2}) no such $F$ exists. This will imply Theorem \ref{keli1} immediately.

To give an interpretation of Theorem \ref{keli1} in the language of graph theory we define the fundamental concept of the paper: A polybox code. 

A set $S$ of arbitrary objects will be called an \textit{alphabet}, and the elements of $S$ will be called \textit{letters}. A permutation $s\mapsto s'$ of the alphabet $S$ such that $s''=(s')'=s$ and $s'\neq  s$ is said to be a \textit{complementation}. 
Each sequence of letters $s_1\ldots s_d$ from the set $S$ is called a \textit{word}. The set of all words of length $d$ is denoted by $S^d$. 
Two words $u=u_1\ldots u_d$ and $v=v_1\ldots v_d$ are \textit{dichotomous} if there is $j\in[d]$ such that $u'_j=v_j$. If $V\subseteq S^d$ consists of pairwise dichotomous words, then we call it a \textit{polybox code} (or \textit{polybox genome}). 
Two words $u,v\in S^d$ form \textit{a twin pair } if there is $j\in[d]$ such that $u'_j=v_j$ and  $u_i=v_i$ for every $i\in [d]\setminus\{j\}$.

If $S$ is an alphabet with a complementation, then {\it a d-dimensional Keller graph on the set $S^d$} is the graph in which two vertices $u,v\in S^d$ are adjacent if they are dichotomous but do not form a twin pair. It is a generalization of the well known {\it d-dimensional Keller graph} in which $S=\{0,1,2,3\}$ and the complementation is given by $0'=2$ and $1'=3$ (\cite{CS2}). 

As we shall show Theorem \ref{keli1} for Keller graphs on $S^d$ reads as follows: 
\begin{tw}
\label{clic}
Let $V\subset S^7$ be a clique in the Keller graph on $S^7$ such that there is $i\in[7]$ and there are four words $v^a,v^b,v^c,v^d\in V$ with $v^a_i\in \{a,a'\},...,v^d_i\in \{d,d'\}$, where $\{a,b,c,d\}\subset S$.  Then $|V|<2^7$. 
\end{tw}

We also show that the mentioned above Theorem \ref{ke2} together with Theorem 29 of \cite{Kis22} give us a clue about the form of a counterexample in the last case of the conjecture, if it exists:
 
\begin{wn}
\label{wk}
There is a counterexample to Keller's conjecture in dimension seven if and only if there is a clique $V$ in the Keller graph on $\{a,a',b,b',c,c'\}^7$ having 128 words such that there are $i\in[7]$ and three words $v^a,v^b,v^c\in V$ with $v^a_i\in \{a,a'\},v^b_i\in \{b,b'\}$, $v^c_i\in \{c,c'\}$ and $|V^{j,l}|>16$ for every $l\in  \{a,a',b,b',c,c'\}$ and every $j\in [7]$ such that $V^{j,l}\neq\emptyset$, where $V^{j,l}=\{v\in V\colon v_j=l\}$.
\end{wn}

At the end of the paper we shall also give a new proof of Keller's conjecture in dimensions $d\leq 6$.

\smallskip
The presented paper is a continuation of two earlier papers \cite{Kis,Kis22} devoted to Keller's conjecture in dimension seven. Therefore, the first two sections containing the basic concepts have been limited to a necessary minimum. More comprehensive presentation of the notions dealing with the structure of polybox codes can be found in the mentioned two papers.

The outline of the paper will be presented after the following section: 

\vspace{-5mm}
\section{Basic notions}

In this section we present the basic notions on dichotomous boxes and words (details can be found in \cite{Kis,Kis22,KP}). We start with systems of boxes.

In the whole paper, if $\ka X$ is a family of sets, then $\bigcup \ka X=\bigcup_{A\in \ka X}A$. Moreover, if $Y$ is a set, then a {\it partition of} $Y$ is a family $\ka Y$ of its pairwise disjoint subsets such that $\bigcup \ka Y=Y$.  

\vspace{-4mm}
\subsection{Dichotomous boxes and  polyboxes}

Let $X_1,\ldots ,X_d$ be non-empty sets with $|X_i|\geq 2$ for every $i\in [d]$. The set $X=X_1\times\cdots \times X_d$ is called a $d$-{\it box}.
A non-empty set $K \subseteq X$ is called a \textit{ box} if $K=K_1\times\cdots \times K_d$ and
$K_i\subseteq X_i$ for each $i\in [d]$.

\smallskip 
The box  $K$ is said to be \textit{ proper} if $K_i\neq X_i$ for each $i\in [d]$.  
Two boxes $K$ and $G$ in $X$ are called \textit{dichotomous}
if there is $i\in [d]$ such that $K_i=X_i\setminus G_i$. A \textit{suit} is any collection of pairwise
dichotomous boxes. A suit is \textit{proper} if it consists of proper boxes. 
A non-empty set $F\subseteq X$ is said to be a \textit{ polybox} if
there is a suit $\ka F$ for $F$, that is, if $\bigcup \ka F=F$. In other words, $F$ is a polybox if it has a partition into pairwise dichotomous boxes. 
A polybox $F$ is {\it rigid} if it has exactly one suit, that is, if $\ka F$ and $\ka G$ are suits for a rigid polybox, then $\ka F=\ka G$. 

 
 A proper suit for a $d$-box $X$ is called a {\it minimal partition } of $X$. In \cite{GKP} we showed that a suit  $\ka F$ is a minimal partition of a $d$-box $X$ if and only if $|\ka F|=2^d$.
 


Two boxes $K,G\subset X$ are said to be a {\it twin pair} if $K_j=X_j\setminus G_j$ for some $j\in [d]$ and $K_i=G_i$ for every $i\in [d]\setminus\{j\}$. Alternatively, two dichotomous boxes $K,G$ are a twin pair if $K\cup G$ is a box. Observe that the suit for a rigid polybox cannot contain a twin pair.

The next concept is of  particular importance in an analysis of the structure of suits.
Let $X$ be a $d$-box, and let $l_i=\{x_1\}\times \cdots \times \{x_{i-1}\}\times X_i \times \{x_{i+1}\}\times \cdots \times \{x_d\}$, where $x_j\in X_j$ for $j\in [d]\setminus\{i\}$. 
A set $F\subseteq X$ is called an {\it $i$-cylinder} (Figure 2) if for every set $l_i$ 
one has 
$$
l_i\cap F=l_i \; \; \; {\rm or}\; \; \; l_i\cap F=\emptyset.
$$

\vspace{-0mm}
{\center
\includegraphics[width=7cm]{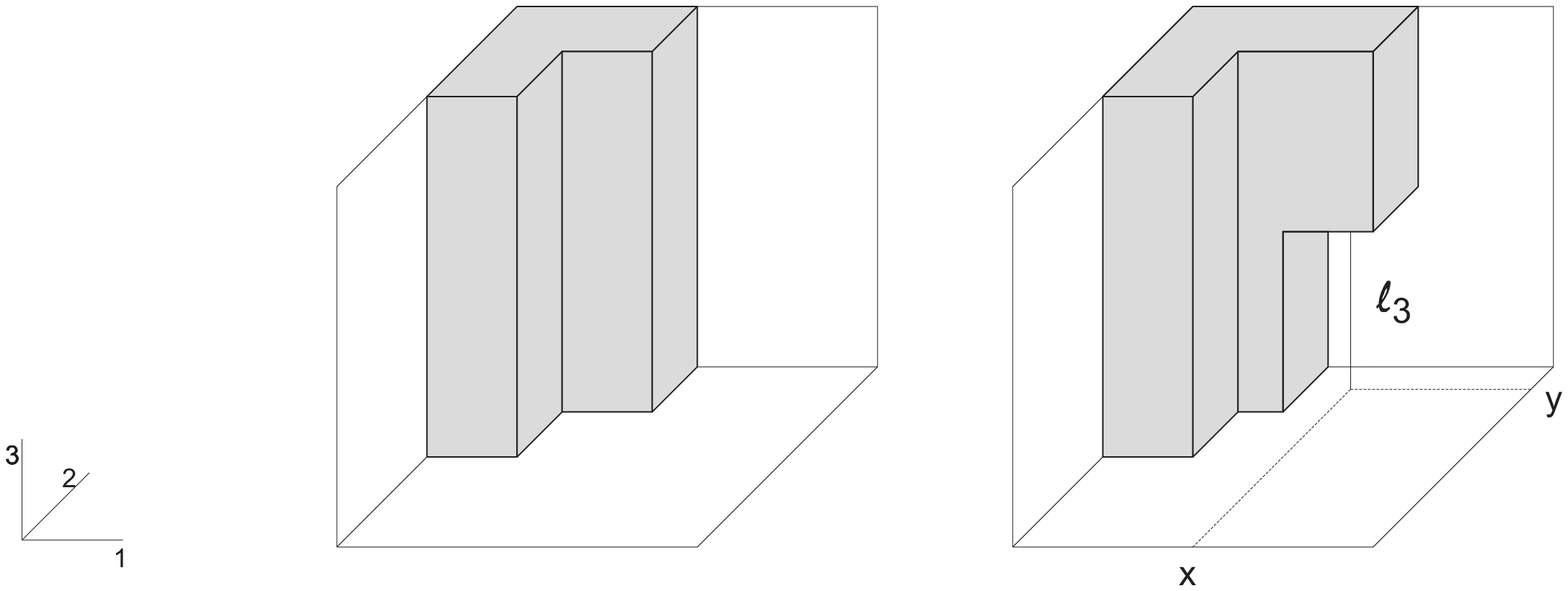}\\
}

\medskip
\noindent{\footnotesize Figure 2: The set on the left is a $3$-cylinder in $X=[0,1]^3$, and the set on the right is not, because the set $l_3=\{x\}\times \{y\}\times [0,1]$ has a non-empty intersection with this set but $l_3$ is not entirely contained in it. 
}

\smallskip
Let $\ka F$ be a suit for a $d$-box $X$, $A\subseteq X_i$, and let $\ka F^{i,A}=\{K\in \ka F\colon K_i=A\}$. Note now that the partition $\ka F$ has the ``cylindrical" structure (\cite{CS1,LS2,La}) (compare Example \ref{p}): Let $i\in [d]$, and let $A^1,...,A^{k_i}\subset X_i$ be all sets such that $\ka F^{i,A^{j}}\neq\emptyset$ for $j\in [k_i]$ and $A^n\not \in \{A^j,(A^j)^c\}$ for $n,j \in [k_i], n\neq j$, where  $(A^j)^c=X_i\setminus A^j$. Then the set $\bigcup (\ka F^{i,A^{j}}\cup \ka F^{i,(A^{j})^c})$ is a non-empty $i$-cylinder in $X$ and $\ka F=\bigcup_{j\in [k_i]} (\ka F^{i,A^{j}}\cup \ka F^{i,(A^{j})^c})\cup \ka F^{i,X_i}$.

\subsection{Cube tilings and dichotomous boxes}

Every two cubes $[0,1)^d+t$ and $[0,1)^d+p$ in an arbitrary cube tiling $[0,1)^d + T$ of $\er^d$ satisfy {\it Keller's condition}: There is $i\in [d]$ such that $t_i-p_i\in \zet\setminus\{0\}$, where $t_i$ and $p_i$ are $i$th coordinates of the vectors $t$ and $p$ (\cite{Ke1}). For any cube $[0,1]^d+x$, where $x=(x_1,...,x_d)\in \er^d$, the family $\ka F_x=\{([0,1)^d+t)\cap ([0,1]^d+x)\neq\emptyset:t\in T\}$ is a partition of the cube $[0,1]^d+x$. Every two boxes  $K,G\in \ka F_x$ are, by Keller's condition, dichotomous: There is $i\in [d]$ such that $K_i$ and $G_i$ are disjoint and $K_i\cup G_i=[0,1]+x_i$. Moreover, since cubes in cube tilings are half-open, every box $K\in \ka F_x$ is proper, and consequently the family $\ka F_x$ is a minimal partition. The structure of the partition $\ka F_x$ reflects the local structure of the cube tiling $[0,1)^d+T$. Obviously, a cube tiling $[0,1)^d+T$ contains a twin pair if and only if the partition $\ka F_x$ contains a twin pair for some $x\in \er^d$ (\cite{La,P}) (see Figure 1). 
Observe also that if $\ka F_x=\ka F_x^{i,A^1}\cup \ka F_x^{i,(A^1)^c}\cup \cdots \cup \ka F_x^{i,A^{k_i(x)}}\cup \ka F_x^{i,(A^{k_i(x)})^c}$, then $|L(T,x,i)|=k_i(x)$ (compare (\ref{re})).

\medskip

\vspace{-5mm}
\subsection{Distribution of words in codes. Realizations of codes}


In the presented paper suits will be encoded with polybox codes described in Introduction.
In what follows we assume that $S$ is a finite alphabet with a fixed complementation.   
If $V\subseteq S^d$, $l\in S$ and $i\in [d]$, then $V^{i,l}=\{v\in V\colon v_i=l\}$. If $S=\{a_1,a'_1,...,a_k,a_k'\}$, then the representation $V=V^{i,a_1}\cup V^{i,a_1'}\cup \cdots \cup V^{i,a_{k}}\cup V^{i,a_{k}'}$
will be called a {\it distribution of words in} $V$.  
Moreover, let $D_i(V)=((|V^{i,a_1}|,|V^{i,a'_1}|),\ldots, (|V^{i,a_k}|,|V^{i,a'_k}|))$ for $i\in [d]$. 
For example, if $S=\{a,a',b,b'\}$ and $V=\{bba,\; baa',\; b'ab,\;b'a'a\}$,
then 

\vspace{-3mm}
$$
D_1(V)=((0,0),(2,2)),\;\; D_2(V)=((2,1),(1,0)),\;\; D_3(V)=((2,1),(1,0)).
$$

A natural interpretation of a polybox code is a suit for a polybox: Let $X=X_1\times \cdots \times X_d$ be a $d$-box.
Suppose that for each $i\in[d]$ a mapping $f_i\colon S\to 2^{X_i}\setminus \{\emptyset,X_i\}$ is such that $f_i(s')=X_i\setminus f_i(s)$. We define the mapping $f\colon S^d\to 2^X$ by $f(s_1\ldots s_d)=f_1(s_1)\times\cdots\times f_d(s_d).$ 
If now $V\subset S^d$ is a code, then the set of boxes $f(V)=\{f(v)\colon v\in V\}$ is a suit for the polybox $\bigcup f(V)$. The set $f(V)$ is said to be a \textit{realization} of the set $V$. A code has infinitely many realizations. But to effectively study the structure of suits we shall use the following realization which has particular nice properties. Let
$$
 ES=\{B\subset S\colon |\{s,s'\}\cap B|=1, \text{whenever $s\in S$}\}\;\;\;{\rm and}\;\;\; E s=\{B\in ES\colon s\in B\}.
$$
Let $V\subseteq S^{d}$ be a polybox code, and let $v\in V$. The {\it equicomplementary realization} of the word $v$ is the box   
$$
\breve{v}=Ev_1\times \cdots \times Ev_d
$$
in the $d$-box $(ES)^d =ES\times \cdots \times ES.$ The equicomplementary realization of the code $V$ is the family
$$
E(V)=\{\breve{v}:v\in V\}.
$$ 
If $s_1,\ldots , s_n\in S$ and $s_i\not\in\{s_j, s'_j\}$ for every $i\neq j$, then
\begin{equation}
\label{dkostki}
|E s_1 \cap\dots\cap E s_n|=(1/2^{n})|ES|.
\end{equation}
The value of the realization $E(V)$, where $V\subseteq S^d$,  lies in the equality (\ref{dkostki}). 
In particular, boxes in $E(V)$ are of the same size: $|E v_i|=(1/2)|ES|$ for every $i\in [d]$ and consequently  $|\breve{v}|=(1/2^d)|ES|^d$ for $v\in E(V)$. Thus, two boxes $\breve v, \breve w\subset (ES)^d$ are dichotomous if and only if $\breve v\cap \breve w=\emptyset.$ The same is true for cubes in a cube tiling of a polycube $F\subset \te^d$ and therefore working with the boxes $\breve{v}, v\in V$,  we can think of them as translates of the unit cube in $\te^d$. 

\vspace{-3mm}
\subsection{Cover of a code, equivalent and rigid polybox codes}
Let $V,W\subseteq S^d$ be polybox codes, and let $w\in S^d$. We say that $w$ is \textit{covered} by $V$, and write $w\sqsubseteq V$, if $\breve{w}\subseteq \bigcup E(V)$. 
If  $w\sqsubseteq V$ for every $w\in W$, then we write $W \sqsubseteq V$, and the code $V$ is called a {\it cover} of the code $W$. Every cover $V$ of $w$ such that $w\not\in V$ has the following useful property (\cite[Theorem 10.6]{KP}): {\it If $v\in V$ is a word such that $\breve{v}\cap \breve{w}\neq\emptyset$ and $|\{i\in [d]\colon v_i\neq w_i\}|\geq |\{i\in [d]\colon u_i\neq w_i\}|$ for every $u\in V$ with $\breve{u}\cap \breve{w}\neq\emptyset$, then there is $p\in V$ such that $p_i\in \{v_i,v_i'\}$ for every $i\in [d]$, $\breve{p}\cap \breve{w}\neq\emptyset$, and moreover the number $|\{i\in [d]\colon v_i= p'_i\}|$ is odd.} This property simplifies the computations of covers of $w\in S^d$ (compare algorithm {\small {\sc CoverCode}} at the end of this section).

Another consequence of (\ref{dkostki}) dealing with covers of words is given in the following lemma (\cite[Lemma 2]{Kis22}):
\begin{lemat}
\label{twin}
Let $V\subset S^d$ be a cover of a word $w\in S^d$ such that $\breve{w}\cap \breve{v}\neq\emptyset$ for $v\in V$. If $V$ does not contain twin pairs, then also the suit $\ka F=\{\breve{w}\cap \breve{v}\colon v\in V\}$ for $\breve{w}$ does not contain such pairs.

\hfill{$\square$}
\end{lemat}

Polybox codes $V,W\subseteq S^d$ are said to be {\it equivalent} if $V \sqsubseteq W$ and  $W \sqsubseteq V$. Thus, $V$ and $W$ are equivalent if and only if  $\bigcup E(V)=\bigcup E(W)$. Obviously, if $V$ and $W$ are equivalent, then $|V|=|W|$. Two codes $V$ and $W$ are {\it disjoint } if $V\cap W=\emptyset$.

Let $q\in S^d$. We say that codes $V,W\subset S^d$ are $q$-{\it equivalent} if 
\begin{equation}
\label{qeq}
\breve{q}\cap \bigcup E(V)=\breve{q}\cap \bigcup E(W).
\end{equation}

In the paper $q$-equivalent codes appear in the strictly defined circumstances, which are explained at the beginning of the next subsection (see Lemma \ref{wac}).  

A polybox code $V\subset S^d$ is called {\it rigid} if there is no code $W\subset S^d$ which is equivalent to $V$ and $V\neq W$. Thus, if polybox codes $V,W$ are equivalent and one of them is rigid, then $V=W$. It can be checked that the code $V$ given in Example \ref{p} is rigid, while $V=\{aaa,a'aa\}$ is not rigid as $W=\{baa,b'aa\}$ is equivalent to $V$.
Observe that, rigid polybox codes cannot contain a twin pair.

\vspace{-5mm}
\subsection{Geometry of dichotomous boxes}

In this subsection we describe the main techniques used in the paper.  

Throughout the paper we consider two disjoint equivalent codes $V,W\subset S^d$. 
If $A=\{i_1<...<i_k\}\subset [d]$ and $v\in S^d$, then $v_A=v_{i_1}\ldots v_{i_k}\in S^k$ and $V_A=\{v_A\colon v\in V\}$. To simplify notation we let $v_{i^c}=v_{\{i\}^c}$, that is, the word $v_{i^c}\in S^{d-1}$ arises from $v$ by skipping the letter $v_i$ in $v$. Moreover, $V_{i^c}=\{v_{i^c}:v\in V\}$.

\smallskip
$\bullet$ {\it The structure of $V$ from the suit for $\breve{w}$, where $w\sqsubseteq V$}. Let $w\sqsubseteq V$. Then $\breve{w}\subseteq \bigcup E(V)$ and the set of boxes $\ka F=\{\breve{w}\cap \breve{v}:v\in V\}$ is a suit for $\breve{w}$. Assume that $w_i=b$ and the sets

\vspace{-3mm}
$$
U^{i,a}=\{v\in V^{i,a}\colon \breve{v}\cap \breve{w}\neq\emptyset\},\;\; U^{i,a'}=\{v\in V^{i,a'}\colon \breve{v}\cap \breve{w}\neq\emptyset\}.
$$
are non-empty (actually, if one of these sets is non-empty, then so is the second one).
The set of boxes $\ka F$ is a suit for the box $\breve{w}$,  which means that the set 

\begin{equation}
\label{walle}
\bigcup (\{\breve{v}\cap \breve{w}: v\in U^{i,a}\} \cup \{\breve{v}\cap \breve{w}: v\in U^{i,a'}\})
\end{equation}
is an $i$-cylinder in the box $\breve{w}$ (compare Lemma 12 of \cite{Kis22} and Example \ref{p}). Therefore, we have
\begin{lemat}
\label{wac}
Let $V,W\subset S^d$ be equivalent polybox codes. If there are $i\in [d]$, $l\in S$ and $w\in W\setminus (W^{i,l}\cup W^{i,l'})$ such that the set $U^{i,l}=\{v\in V^{i,l}\colon \breve{v}\cap \breve{w}\neq\emptyset\}$ is non-empty, then the codes  $U^{i,l}_{i^c}, U^{i,l'}_{i^c}$ are $w_{i^c}$-equivalent. \hfill{$\square$}
\end{lemat}


This is mentioned above context in which $q$-equivalent codes will appear in the paper. It is worth analyzing the following example in which we describe  the typical situation encountered in the paper.

\vspace{-2mm}
\begin{ex}
\label{p}
{\rm In Figure 3 the five boxes on the left are a realization of the polybox code $V=\{aaa,a'a'a',baa',a'ba, aa'b\}$, and the box in the middle is a realization of the word $w=bbb$. Since the set $V$ is a cover of $w$, that is, $w\sqsubseteq V$, we have $\breve{w}\subset \bigcup E(V)$. Thus, the $3$-box $\breve{w}$ is divided into pairwise dichotomous boxes of the form $\breve{v}\cap \breve{w}$ for $v\in V$. (In other words, the set of boxes $\ka F=\{\breve{v}\cap \breve{w}:v\in V\}$ is a suit for $\breve{w}$.)  The set  $\bigcup (\{\breve{v}\cap \breve{w}: v\in U^{3,a}\} \cup \{\breve{v}\cap \breve{w}: v\in U^{3,a'}\})$,
 where  $U^{3,l}=\{v\in V^{3,l}: \breve{v}\cap \breve{w}\neq\emptyset\}$ for $l\in \{a,a'\}$,
 is a $3$-cylinder in the box $\breve{w}$. Therefore, $\bigcup \{\breve{v}_{3^c}\cap \breve{w}_{3^c}: v\in U^{3,a}\}=\bigcup \{\breve{v}_{3^c}\cap \breve{w}_{3^c}: v\in U^{3,a'}\}$, and thus the codes $U^{3,a}_{3^c},U^{3,a'}_{3^c}$ are $w_{3^c}$-equivalent.}

\end{ex}

{\center
\includegraphics[width=10cm]{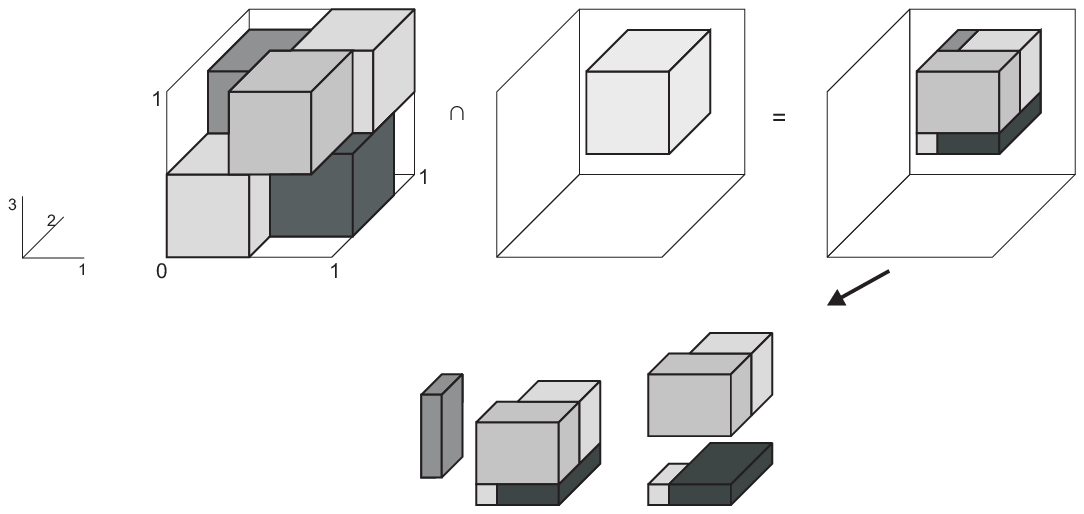}\\
}

\medskip
\noindent{\footnotesize Figure 3: On the top: The light box (the box in the middle, which may be interpreted as $\breve{w}$) is contained in the sum of five pairwise dichotomous boxes (the boxes on the left). These boxes determine a partition of the light box into pairwise dichotomous boxes (the partition on the right). On the bottom: The boxes in this partition are arranged into $3$-cylinders. The sets $Ea$ and $Eb$ are identified with $[0,\frac{1}{2})$ and $[\frac{1}{4},\frac{3}{4}]$, respectively.
  }

\medskip
$\bullet$  {\it The  structure of $W$ from the distribution of words in $V$}.  Below, in (\textbf{P}), (\textbf{C}) and (\textbf{Co}) we show how to use an information on a  distribution of words in $V$ to say something about the distribution of words in $W$.

\vspace{0mm}
{\center
\includegraphics[width=11cm]{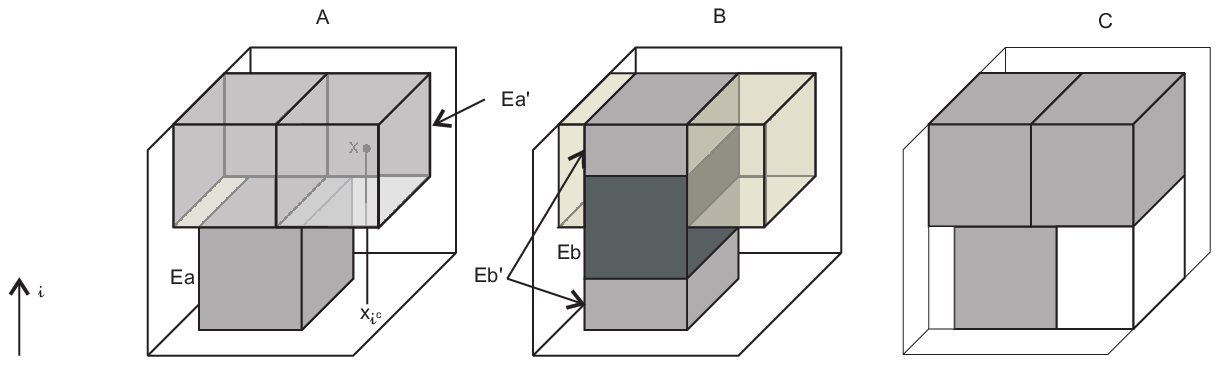}\\
}

\medskip
\noindent{\footnotesize 
Figure 4: A scheme of realizations $E(V)$ (A), $E(W)$ (B), where $V=V^{i,a}\cup V^{i,a'}$ and $W=W^{i,a'}\cup W^{i,b}\cup W^{i,b'}$. We assume that the codes $V$ and $W$ are equivalent, and thus $\bigcup E(V)=\bigcup E(W)$.    
  }
  
\medskip

Let $V,W\subseteq S^d$ be polybox codes and assume that $V$ and $W$ are equivalent. Recall that $V^{i,l}_{i^c}=(V^{i,l})_{i^c}$ for every $i\in [d]$ and $l\in S$. Moreover, $x_{i^c}=(x_1,...,x_{i-1},x_{i+1},...,x_d)$ for $x\in (ES)^d$.

\medskip
\noindent
(\textbf{P}): {\bf Projections}. Suppose that there is $x\in \bigcup E(V^{i,l'})$ such that $x_{i^c}\not\in \bigcup E(V^{i,l}_{i^c})$ (see Figure 4A, where $l=a$). Since $\bigcup E(V)=\bigcup E(W)$ and (\ref{dkostki}), the point $x$ can be covered only by a box $\breve{w}\in E(W)$ such that $w_i=l'$. Thus, $W^{i,l'}\neq\emptyset$. In particular, if $W^{i,l}=\emptyset$ and $W^{i,l'}\neq\emptyset$, then $W^{i,l'}\sqsubseteq V^{i,l'}$ (compare Figures 4A and 4B). Note also that, by (\ref{dkostki}), the white box in Figure 4C cannot be intersected by a box $\breve{v}$ for $v\in V\setminus (V^{i,a}\cup V^{i,a'})$. Thus, the white box cannot be intersected by a box $\breve{w}$ for $w\in W$. Finally, observe that if $\breve{w}\cap \breve{v}\neq\emptyset$ for some $w\in W$ and $v\in V$ with $w_i\not\in \{v_i,v_i'\}$, then every point $x\in \breve{v}\setminus \breve{w}$ such that $x_{i^c}\in \breve{w}_{i^c}$ is covered only by boxes $\breve{u}$, $u\in W$, such that $u_i=w_i'$.

\medskip
It follows from the above

\begin{lemat}
\label{rr1}
Let $V,W\subset S^d$ be equivalent codes. Then 
\begin{equation}
\label{dif}
\bigcup E(V^{i,l}_{i^c})\setminus \bigcup E(V^{i,l'}_{i^c})=\bigcup E(W^{i,l}_{i^c})\setminus \bigcup E(W^{i,l'}_{i^c})
\end{equation}
for every $i\in [d]$ and every $l\in S$.\hfill{$\square$}
\end{lemat}

For example, it can be computed that if $V,W\subset S^d$, $d=5,6$, are disjoint equivalent codes without twin pairs and $|W^{1,l}|=5$ for $l\in \{a,a'\}$ and $|V^{1,l}|=1$ for $l\in \{a,a'\}$, then, up to isomorphism (it will be defined in Subsection 2.6), there is only one pair $(W^{1,a},W^{1,a'})$ and $(V^{1,a},V^{1,a'})$ for which (\ref{dif}) is valid:
$$
W^{1,a}=\{aaaabb, aa'a'a'bb, abaa'bb, aa'babb, aaa'bbb\},\;\; V^{1,a}=\{abbbbb\},
$$
\begin{equation}
\label{5n5}
W^{1,a'}=\{a'ab'abb, a'bbb'bb, a'bb'a'bb, a'b'a'a'bb, a'b'babbb\}, \;\;V^{1,a'}=\{a'aa'b'bb\}.
\end{equation}
Note that $v\not\sqsubseteq W^{1,l}$ for $v \in W^{1,l'}$ and $l\in \{a,a'\}$. (We shall use this fact in Section 6.)

\smallskip
Let $V,W\subset S^d$ be codes, and let 
$$
m(V\setminus W)=\left(\frac{2}{|ES|}\right)^{d}(|\bigcup E(V)\setminus\bigcup E(W)|).
$$
The number $m(V\setminus W)$ counts the fraction of boxes $\breve{v},v\in V$, which are not covered by boxes $\breve{w},w\in W$.
It follows from Lemma \ref{rr1} that for every equivalent codes $V$ and $W$ we have
\begin{equation}
\label{mb}
m(V^{i,l}_{i^c}\setminus V_{i^c}^{i,l'})=m(W^{i,l}_{i^c}\setminus W_{i^c}^{i,l'})\;\; {\rm for}\;\; i\in [d]\;\; {\rm and} \;\; l\in S.
\end{equation}
For example, for the code $V$ in Example \ref{p} we have $m(V^{3,l}_{3^c}\setminus V^{3,l'}_{3^c})=\frac{3}{4}$ for $l\in \{a,a'\}$ (recall that $V^{i,l}_{i^c}\subset S^{d-1}$ for $V\subset S^d$.)

\smallskip
\noindent
(\textbf{C}): {\bf Cylinders}. Suppose that $V^{i,l}\cup V^{i,l'}=\emptyset$ and $W^{i,l}\cup W^{i,l'}\neq\emptyset$ for some $l\in S$. Then, by (\ref{dif}), $\bigcup E(W^{i,l}_{i^c})=\bigcup E(W^{i,l'}_{i^c})$, and hence the set $\bigcup E(W^{i,l}\cup W^{i,l'})$ in an $i$-cylinder in the $d$-box $(ES)^d$, that is, $\bigcup E(W^{i,l}_{i^c})=\bigcup E(W^{i,l'}_{i^c})$ (compare Figures 4A and 4B, where $l=b$). 
Thus, by the definition,  the codes $W^{i,l}_{i^c}$ and $W^{i,l'}_{i^c}$ are equivalent.


\smallskip
\noindent
(\textbf{Co}): {\bf Covers}. Suppose that polybox codes $V^{i,l}_{i^c}$ and $W^{i,l}_{i^c}\cup W^{i,s_1}_{i^c}\cup \cdots \cup W_{i^c}^{i,s_k}$ are equivalent, where $s_n\not\in \{l,l',s_j,s_j'\}$ for every $n,j\in [k],n\neq j$. Then 
$$
W^{i,s_1}_{i^c}\cup \cdots \cup W_{i^c}^{i,s_k}\sqsubseteq V^{i,l'}_{i^c}\quad   {\rm and}\quad  V^{i,l}_{i^c}\sqsubseteq  W^{i,l}_{i^c}\cup W^{i,s_1'}_{i^c}\cup \cdots \cup W^{i,s_k'}_{i^c}.
$$
(In Figures 4A and 4B the codes  $V^{i,a}_{i^c}$ and $W^{i,a}_{i^c}\cup W^{i,b}_{i^c}$ are equivalent, where $a=l,b=s_1$, and $W^{i,a}=\emptyset$). Indeed, since boxes in $E(V)$ are pairwise dichotomous and, by the definition of equivalent codes, $\breve{w}_{i^c}\subseteq \bigcup E(V^{i,l}_{i^c})$ for every $w\in W^{i,s_1}\cup \cdots \cup W^{i,s_k}$,  it follows that each point $x\in \breve{w}\setminus \bigcup E(V^{i,l})$ has to be covered by the set $\bigcup E(V^{i,l'})$. Therefore,  $\breve{w}_{i^c}\subseteq \bigcup E(V^{i,l'}_{i^c})$ for every $w\in W^{i,s_1}\cup \cdots \cup W^{i,s_k}$, and consequently  $w_{i^c}\sqsubseteq V^{i,l'}_{i^c}$ for $w\in W^{i,s_1}\cup \cdots \cup W^{i,s_k}$. Thus, $W^{i,s_1}_{i^c}\cup \cdots \cup W_{i^c}^{i,s_k}\sqsubseteq V^{i,l'}_{i^c}$. In the same manner we show that  $V^{i,l}_{i^c}\sqsubseteq  W^{i,l}_{i^c}\cup W^{i,s_1'}_{i^c}\cup \cdots \cup W^{i,s_k'}_{i^c}$.

\medskip
$\bullet$  {\it The structure of $W$ from the structure of $V^{i,l}\cup V^{i,l'}$}.
The following lemma describes a useful relationship between codes $V^{i,l}\cup V^{i,l'}$ and $W^{i,s}\cup W^{i,s'}$, where $s\not\in \{l,l'\}$.

\begin{lemat}
\label{ci}
Let $V,W\subset S^d$ be disjoint equivalent polybox codes, and let $U^{i,l}=\{v\in V^{i,l}\colon \breve{v}\cap \bigcup E(W\setminus (W^{i,l}\cup W^{i,l'}))\neq\emptyset\}$ be non-empty set for some $i\in [d]$ and $l\in S$. Assume that there are disjoint sets $A,B$ such that $[d]\setminus \{i\}=A\cup B$ and there are two words $p,q$ such that:

\smallskip
(i)\;\;\; $w_A=p$ for every $w\in W\setminus (W^{i,l}\cup W^{i,l'})$ with $\breve{w}\cap \bigcup E(U^{i,l}\cup U^{i,l'})\neq\emptyset$.

\smallskip
(ii)\;\;\; $(U^{i,l}\cup U^{i,l'})_B=\{q\}$.

\smallskip
If $P=\{w\in  W\setminus (W^{i,l}\cup W^{i,l'})\colon \breve{w}\cap \bigcup E(U^{i,l}\cup U^{i,l'})\neq\emptyset\}$,
then for every $s\in S\setminus \{l,l'\}$ the codes $P^{i,s}_{i^c}$ and $P^{i,s'}_{i^c}$ are $r$-equivalent, where $r_A=p$ and $r_B=q$. 
In particular, if $|B|\leq 1$, then there is a twin pair in $V$ or in $W$. 
\end{lemat} 
\proof For simplicity, let $i=1$, $l=a$, $s=b$ and $A=\{2,...,k\}, B=\{k+1,...,d\}$. Moreover, $x_B=(x_{k+1},...,x_d)$ and $x_{1^c}=(x_2,...,x_d)$ for $x\in (ES)^d$. To show that that $P^{1,b}_{1^c}$ and $P^{1,b'}_{1^c}$ are $r$-equivalent it is enough to show, by ($i$), that
$$
\bigcup_{w \in P^{1,b}} \breve{w}_{B}\cap \breve{q}=\bigcup_{w \in P^{1,b'}} \breve{w}_{B}\cap \breve{q}.
$$
Take  a point $x_{B}\in \bigcup_{w \in P^{1,b}} \breve{w}_{B}\cap \breve{q}$. Then $x_{B}\in \breve{w}_{B}\cap \breve{q}$ for some $w\in P^{1,b}$. By the definition of $P^{1,b}$, there is $v\in U^{1,a}\cup U^{1,a'}$ such that $\breve{w}\cap \breve{v}\neq\emptyset$. Assume that $v\in U^{1,a}$ and take $y\in \breve{v}\cap \breve{w}$ such that $y_B=x_B$, which is possible by ($ii$). Then, by ({\bf P}), a point $z\in \breve{w}$ such that $z_1\in Ea'\cap Eb$ and $z_{1^c}=y_{1^c}$ must be covered by a box $\breve{u}$ for some $u\in U^{1,a'}$. We take one more point $t\in \breve{u}$ such that $t_1\in Ea'\cap Eb'$ and $t_{1^c}=z_{1^c}$. Clearly, again by ({\bf P}), $t\in \breve{w}^1\cap \breve{u}$ for some $w^1\in P^{1,b'}$. Since $t_B=x_B$ and $\breve{u}_B=\breve{q}$, we obtain $x_B\in \breve{w}^1_B\cap \breve{q}\subseteq \bigcup_{w \in P^{1,b'}} \breve{w}_{B}\cap \breve{q}$. Consequently,  $\bigcup_{w \in P^{1,b}}\breve{w}_{B}\cap \breve{q}\subseteq  \bigcup_{w \in P^{1,b'}} \breve{w}_{B}\cap \breve{q}$. In the same manner we show the reverse inclusion.

To prove the second part of the lemma let first $|B|=1$. If $|P^{1,b}|\geq 2$, then there are two words $w,u\in W^{1,b}$ such that $w_1=u_1=b$ and $w_A=u_A$. Since $w,u$ are dichotomous, $w_d=u_d'$, and thus $w$ and $u$ are a twin pair.  If $P^{1,b}=\{w\}$, $P^{1,b'}=\{u\}$ and $v\in U^{1,a}\cup U^{1,a'}$ is such that $\breve{w}\cap \breve{v}\neq\emptyset$ and $\breve{u}\cap \breve{v}\neq\emptyset$, then, by Lemma \ref{wac} (in which we change the role of $V$ and $W$) the codes $P^{1,b}_{1^c}$ and $P^{1,b'}_{1^c}$ are $v_{1^c}$-equivalent, and hence, by (\ref{dkostki}), $w_d=u_d$. Thus, $w,u$ form a twin pair.

If $B=\emptyset$, then for every $w\in P^{1,b}$ and $u\in P^{1,b'}$, we have $w_A=u_A$ and $w_1=u_1'$, that is, the words $w,u$ form a twin pair.  
\hfill{$\square$}

\medskip
Let  
$g\colon S^d\times S^d\to \zet$ be defined by the formula $g(v,w)=\prod^d_{i=1}(2[v_i=w_i]+[w_i\not\in\{v_i, v'_i\}])$
where $\iver {p}=1$ if the sentence $p$ is true and $\iver {p}=0$ if it is false. Let $w\in S^d$, and let $V\subset S^d$ be a polybox code. In \cite{KP} it was showed that  
\begin{equation}
\label{2d}
w \sqsubseteq V \Leftrightarrow \sum_{v\in V} g(v,w)= 2^d.
\end{equation}

\medskip
$\bullet$ {\it The structure of $V$ from $($\ref{2d}$)$}. Let $V,W\subset S^d$ be disjoint equivalent polybox codes.
Then for every $w\in W$ we have $w\sqsubseteq V, w\not\in V$, and, by (\ref{2d}),  $\sum_{v\in V} g(v,w)= 2^d$, where $g(v,w)\in \{0,1,2,\ldots ,2^{d-1}\}$ for every $v\in V$. Assume that $w=b\ldots b$ and let $\{v^1,\ldots ,v^k\}\subseteq V$ be such that $\breve{w}\cap \breve{v}^i\neq\emptyset$ for every $i\in [k]$ and $w\sqsubseteq \{v^1,\ldots ,v^k\}$. The solutions of the system of equations $\sum_{i=0}^{d-1}x_i2^{i}=2^d, \sum_{i=0}^{d-1}x_i=k$, where $x_i$ are non-negative integers for $i\in \{0,1,\ldots ,d-1\}$, show the frequency of the letter $b$ in the words from the set $\{v^1,\ldots ,v^k\}$. We explain this on the following example. Recall first that $g(v,w)= 2^{i}$ if and only if $v_j=b$ for every $j\in I\subset [d]$, $|I|=i$ and $v_j\not \in \{b,b'\}$ for $j\in [d]\setminus I$. In the example we assume that $d=3$, $w=bbb$ and $k=5$. The above system has two solutions: $x_0=2,x_1=3,x_2=0$ and $x_0=4,x_1=0,x_2=1$. It follows from the first solution that in the cover $\{v^1,\ldots ,v^5\}$ of $w$ there are exactly three words such that each of them contains exactly one letter $b$ and two words which have no letter $b$ or, by the second solution,  in the set $\{v^1,\ldots ,v^5\}$ there is exactly one word with two letters $b$ and the rest four words have no letter $b$. This observation is quite useful in the computations of covers of a word $w\in W$ as it allows us to restrict the number of words which have to be considered during the computations (see algorithm {\small {\sc CoverWord}} at the end of this section).

\medskip
$\bullet$ {\it The structure of $V$ from a graph of siblings}

\smallskip
In \cite{Kis} we defined a graph on a polybox code $V$. We now recall the definition of it.

Two words $v,u\in S^d$ such that $v_i\not\in \{u_i,u_i'\}$ for some $i\in [d]$, and  $v_{i^c}$, $u_{i^c}$ is a twin pair are called $i$-{\it siblings} (the two top boxes in Figure 3 (on the left) are a realization of $1$-siblings). 

Let $V\subseteq S^d$ be a polybox code.  {\it A graph of siblings on V} is a graph $G=(V,\ka E)$ 
in which two vertices $v,u\in V$ are adjacent if they are $i$-siblings for some $i\in [d]$. We colour each edge in $\ka E$ with the colours from the set $[d]$: An edge $(v,u)\in \ka E$ has a colour $i\in [d]$ if $v,u$ are $i$-siblings. 
Observe that if $v,u$ are $i$-siblings in a polybox code $V$ such that $v_i=l$ and  $u_i=s$, $l\not\in \{s,s'\}$, then the set $V^{i,l}_{i^c}\cup V^{i,s}_{i^c}$ contains the twin pair $v_{i^c},u_{i^c}$. In Section 5 we shall use a subgraph of $G$ to estimate the cardinalities of covers of $w\in W$ by words from $V$.

\subsection{Isomorphic polybox codes}

   
If $v\in S^d$, and $\sigma$ is a permutation of the set $[d]$, then $\sigma^*(v)=v_{\sigma(1)}\ldots v_{\sigma(d)}$. For every $i\in [d]$ let $h_i:S\rightarrow S$ be a bijection such that $h_i(l')=(h_i(l))'$ for every $l\in S$, and let $h:S^d\rightarrow S^d$ be defined by the formula $h(v)=h_1(v_1)\ldots h_d(v_d)$. We say that polybox codes $P,Q\subset S^d$ are {\it isomorphic} if there are $\sigma$ and $h$ such that $Q=\{h_1(v_{\sigma(1)})\ldots h_d(v_{\sigma(d)}): v\in P\}$. The composition $h\circ \sigma^*$ is an {\it isomorphism} between $P$ and $Q$. Let $V$ and $W$ be disjoint, equivalent and twin pair free polybox codes, and let $h\circ \sigma^*(V)$ be an isomorphic code to $V$. It follows from the definition of the isomorphism $h\circ \sigma^*$ that the codes $h\circ \sigma^*(V)$ and  $h\circ \sigma^*(W)$ are also disjoint, equivalent and do not contain a twin pair. To show this, it is enough to notice that the definition of $h\circ \sigma^*$ guarantees that  $h\circ \sigma^*(w)\sqsubseteq h\circ \sigma^*(V)$, whenever $w\sqsubseteq V$. Moreover, for every $i\in [d]$ and $l\in S$ there are $j\in [d]$ and $s\in S$ such that $|V^{i,l}|=|(h\circ \sigma^*(V))^{j,s}|$.
Therefore, usually we shall assume that $V$ or $W$ contains specific codes. In the following example we explain this statement presenting a typical situation that appears in the paper:

\vspace{-4mm}
\begin{ex}
\label{izop}
{\rm Let $V,W$ be two disjoint twin pair free equivalent codes. Suppose that $V$ contains a code $Q$ which is isomorphic to a code $U=\{aabbbb,aa'abbb,a'babbb,a'aa'bbb\}$. If $f$ is a isomorphism between $Q$ and $U$, then $U\subset f(V)$. Assume now that from this inclusion we are able to deduce that $f(W)$ has to contain a code $R$ such that: $R=R^{1,b}\cup R^{1,b'}$, $R_A=\{bb\}$ and the codes $R^{1,b}_B$, $R^{1,b'}_B$ are $q$-equivalent for $q=bbb$, where $A=\{2,3\}$ and $B=\{4,5,6\}$. Let $K$ and $M$ be $q$-equivalent codes such that there is an isomorphism $h$ with $h_{|_B}(q)=q$, $h_{|_B}(R^{1,b}_B)=K$, $h_{|_B}(R^{1,b'}_B)=M$ and moreover $h_{|_H}={\rm id}$, where $H=\{1,2,3\}$. Then $U\subset h\circ f(V)$ and $h\circ f(W)$ contains a code $T$ such that $T=T^{1,b}\cup T^{1,b'}$, $T_H=R_H$ and $T^{1,b}_B=K$, $T^{1,b'}_B=M$. Clearly, $h\circ f(V)$ and $h\circ f(W)$ are disjoint twin pair free equivalent codes and for every $i\in [d]$ and $l\in S$ there are $j\in [d]$ and $s\in S$ such that $|V^{i,l}|=|(h\circ f(V))^{j,s}|$ and $|W^{i,l}|=|(h\circ f(W))^{n,r}|$ for some $n\in [d]$ and $r\in S$. Therefore, we may assume at the very beginning that $U\subset V$ and $T\subset W$.  

}

\end{ex}

\subsection{Algorithms} 
 
Below we give two algorithms used in the paper. 
 
\medskip
{\bf Algorithm} {\small {\sc CoverWord}}.

\medskip
Let $u\in S^d$, $d=5,6$, be a word, $k\geq 5$ be an integer, and let $\ka C^k$ be the family of all $k$-elements twin pair free covers of $u$ such that $\breve{u}\cap \breve{v}\neq\emptyset$ for every $v\in C$ and every $C\in \ka C^k$. For simplicity  let $u=b...b$. Let $V_{3,0}=\{aaaa...a,a'a'a'a...a\}$, $V_{3,1}=\{aaaba...a,a'a'a'ba...a\}$, $V_{3,2}=\{aaabba...a,a'a'a'bba...a\}$, $V_{3,3}=\{aaabbb, a'a'a'bbb\}$ and $V_{5,0}=\{aaaaaa...a,a'a'a'a'a'a...a\}$, $V_{5,1}=\{aaaaab,a'a'a'a'a'b\}$. By the property given at the beginning of Subsection 2.4, every cover in $\ka C^k$ contains, up to isomorphism, a code $V_{n,i}$ for some $n\in \{3,5\}$, $i\in \{0,...,3\}$. In the algorithm we use also the property discussed below (\ref{2d}). Our goal is to find the family $\ka C^k$.

\smallskip
{\it Input}. The word $b...b\in S^d$ and the number $k$.

{\it Output}. The family $\ka C^k$.

\smallskip
1. Let $\ka S_k=\{(x_0,...,x_{d-1}) \in \en^d\colon  \sum_{i=0}^{d-1}x_i2^i=2^d\; {\rm and}\; \sum_{i=0}^{d-1}x_i=k\}$, where $\en=\{0,1,2,...\}$.

2. For $i \in \{0,...,d-1\}$ indicate the set $\ka A_i$ consisting of all words $v\in S^d$ such that $v$ contains precisely $i$ letters $b$ and $v$ does not contain the letter $b'$. 


3. Fix $x\in \ka S_k$ and let $s(x)=\{i_1<\cdots <i_m\}$ consists of all $i_j\in \{0,...,d-1\}, j\in [m]$, for which $x_{i_j}>0$. Fix $V_{n,i_1}\subset \ka A_{i_1}$. 
For $i\in s(x)$ let $\ka B_i=\{v\in \ka A_i\colon V_{n,i_1}\cup \{v\}\;\; {\rm is\; a \; twin\;\; pair\;\; free\;\; code}\}$.


4. Let $I$ be the multiset containing $i_1$ with the multiplicity $x_{i_1}-2$ (recall that $x_{i_1}\geq 2$) and  $i_j$ with the multiplicity $x_{i_j}$ for $j\in \{2,...,m\}$. By $I[j]$ we denote the $j$th element of $I$. 

5. Let $\ka D^2=\{V_{n,i_1}\}$. 

6. For $l\in \{2,...,k-1\}$ having computed $\ka D^l$ we compute the set $\ka D^{l+1}$: For $v\in \ka B_{I[l]}$ and for $U\in \ka D^l$ if $U\cup \{v\}$ is a twin pair free code, then we attach it to  $\ka D^{l+1}$. 

7. Clearly, $\ka D^k=\ka D^k(V_{n,i_1},x)$ so let $\bar{\ka C}^k$ be the union of the sets $\ka D^k(V_{n,i_1},x)$ over $x\in \ka S^k$ and  $V_{n,i_1}\subset \ka A_{i_1}$ (recall that $V_{n,i_1}$ depends on $x$).

8. $\ka C^k=\bigcup_{f\in F_b}f(\bar{\ka C}^k)$, where $F_b$ consists of all isomorphism $f$, defined in Subsection 2.6, such that $f(b...b)=b...b$.

\medskip
{\bf Algorithm} {\small {\sc CoverCode}}.

\medskip
Let $U=\{u^1,...,u^n\}$ be a code and for every $i\in [n]$ let $P_i$ be a code such that $\breve{u}^i\cap \breve{p}\neq\emptyset$ for every $p\in P_i$. Let $\ka C_{u^i,P_i}$, $i\in [n]$, be the family of all covers $C_{u^i}$ of the word $u^i\in U$ such that $P_i\subset C_{u^i}$.
Our goal is to find the family $\ka C_U$ of all covers $C_U$ of the code $U$ such that $P_i\subset C_{u^i}$ for $i\in [n]$ and $|C_U|\leq m$ for a fixed $m\in \{1,2...\}$, where $C_{u^i}\subset C_U$ is the cover of the word $u^i$.

\medskip   
{\it Input}. The codes $U$, $P_i$, $i\in [n]$, the number $m$ and the family $(\ka C_{u^i,P_i})_{u^i\in U}$.

{\it Output}. The family $\ka C_U$.

\smallskip
1. For $C_1\in \ka C_{u^1,P_1}$ and  $C_2\in \ka C_{u^2,P_2}$ if the set $C_1\cup (C_2\setminus C_1)$ is a twin pair free code (it is obviously a cover of the code $\{u^1,u^2\}$) and has at most $m$ words, then it is attached to the set $\ka C_{1,2}$.  

2. Assuming that the set $\ka C_{1,...,k}$ has already been computed for $2\leq k<n$, we compute the set $\ka C_{1,...,k,k+1}$: for $C \in \ka C_{1,...,k}$ and  $C_{k+1}\in \ka C_{u^{k+1},P_{k+1}}$ if the set $C\cup (C_{k+1}\setminus C)$ has at most $m$ words and it is a twin pair free code (being a cover of $\{u^1,...,u^{k+1}\}$), then it is attached to the set $\ka C_{1,...,k,k+1}$.

3. $\ka C_U=\ka C_{1,...,n}$.

\section{Covers of a code and $q$-equivalent codes}

To show that every cube tiling $[0,1)^7+T$ of $\er^7$ such that $r^+(T)\geq 6$ contains a twin pair, in \cite[Theorem 2.7]{Kis} we proved

\begin{tw}
\label{t12}
If $V,W\subset S^d$ are disjoint equivalent codes without twin pairs, then $|V|\geq 12$.
\end{tw}

To prove the main theorem of the presented paper (Theorem \ref{keli1}, which resolves the case  $r^+(T)=4$) we have to settle whether there are two disjoint twin pair free equivalent codes $V,W\subset S^d$ for $d=6$ such that $|V|\in \{13,...,16\}$. We shall examine such codes $V,W\subset S^d$ for $d=5$ and $d=6$ separately. The reason is that our method of analyzing the structure of equivalent codes $V$,$W$  works better when for every $i\in [d]$ the sets  $V^{i,l}$ and $W^{i,l}$ are non-empty for at least two letters $l\in S$. Therefore, we shall exclude the case when  $V,W$ are {\it flat}, that is, they are of the form $V=V^{i,l}$ and $W=W^{i,l}$ for some $i\in [d]$ and some $l\in S$. (A twin pair free disjoint equivalent codes $V,W\subset S^d$ with $|V|=12$ are not flat only for $d=4$ (see \cite{Kis22}).)


In the next four sections we show which  codes $V,W$   may be excluded from the considerations.
Those reductions will enable us to carry out the main computations, described in Section  7, using home PC equipped with processor Intel i$7$ and  32{\small GB} {\small RAM} memory.

In Section 3 we give covers of small codes $U\subset S^d$  for $d=4,5$ as well as $q$-equivalent codes which will be used throughout the paper.
In the next section we eliminate some codes $V,W$ based on their distributions of words (Lemma \ref{rr2} and \ref{rr3}$i$). Additionally, in Lemma \ref{rr3}$ii$, we show how the distribution of words in $V$ affects the position of some words in $W$ relative to the words in $V$. In Section 5 we estimate the cardinalities of covers of a single word in $V$ and $W$ (Lemma \ref{le12}). To do that, we first give a result on the rigidity of a code (Lemma \ref{rig}). In Section 6 we show that codes $V,W$ under consideration may be written down in the alphabet $S=\{a,a',b,b'\}$. Finally, in Section 7 based on the reductions made in the previous sections we prove Theorem \ref{keli1}. 

\subsection{Covers of a code}
In the first part of the paper we shall need directly only two, up to isomorphism, twin pair free covers $C_v$ of a word $v\in S^d$, $d=4,5$, such that $|C_v|\in \{5,6\}$ and $\breve{v}\cap \breve{u}\neq\emptyset$ for every $u\in C_v$. They can be easily computed using algorithm {\small {\sc CoverWord}} and for $v=bbbbb$ they are of the form:

\vspace{-7mm}
\begin{center}
\begin{displaymath}
\begin{tabular}{|l|}
\hline
  $C^1=\{aaabb,\;a'a'a'bb,\;baa'bb,\;a'babb,\;aa'bbb\}$ \\ \cline{1-1}
  $C^2=\{aaaab,\;a'a'a'ab,\;baa'ab,\;a'baab,\;aa'bab,\;bbba'b\}$  \\ \cline{1-1}
 \hline
\end{tabular}
\end{displaymath}
\end{center}

\noindent{\footnotesize Table 1: Twin pair free covers of the word $bbbbb$ with five and six words.

}

\smallskip
Obviously, $C^1_{5^c},C^2_{5^c}$ are covers of $v=bbbb$.

\begin{lemat}
\label{po7}
Let $U_i,C_i\subset S^d$, where $S=\{a,a',b,b',c,c'\}$, be twin pair free disjoint codes such that $U_i\sqsubseteq C_i$ and $|U_i|=i$ for $i\in [4]$. Then

\smallskip
$(i)$ $|C_1|\geq 5$ for $d\geq 3$. 

\smallskip
$(ii)$ $|C_2|\geq 7$ for $d=5$ and there is precisely one, up to isomorphism,  cover $C_2$ of $U_2$ with $|C_2|=7$. This is $C_2=\{ba'bbb,aabbb,a'aaaa',a'aaba,a'aa'ab,a'aa'a'a,a'aba'a'\}$, where $U_2=\{bbbbb,b'abbb\}$.


\smallskip
$(iii)$ $|C_2|\geq 8$ for $d=4$ and $|C_3|\geq 9$, $|C_4|\geq 10$ for $d=4,5$. 

\end{lemat}
\proof
We sketch the prove of the fact $|C_4|\geq 10$ for $d=5$. For $n\in [5]$ let $\ka P_n$ be a set of all, up to isomorphism, twin pair free two elements codes $\{v,w\}\subset S^5$ such that $|\{i\in [5]\colon v_i=w'_i\}|=n$, and let $\ka P=\bigcup_{n\in [5]} \ka P_n$. 
Using algorithm {\small {\sc CoverWord}} we compute the family $\ka C_v$ consisting of all twin pair free covers $C_v$ of each word $v$ that appears in $\ka P$ such that $\breve{v}\cap \breve{u}\neq\emptyset$ for every $u\in C_v$ and $|C_v|\leq 9$ for every $C_v\in \ka C_v$.  Next, based on algorithm {\small {\sc CoverCode}}, we compute twin pair free covers $C_Q$ of $Q\in \ka P$ which are disjoint with $Q$ and  $|C_Q|\leq 9$. (Note that computing $C_Q$ we can omit all covers $C_v$, $v\in Q$, such that $|C_v|=9$ and $Q\not\sqsubseteq C_v$.)  For $Q\in \ka P_3\cup \ka P_4\cup \ka P_5$ there are no such covers $C_Q$, and every non-empty cover $C_Q$ of $Q\in \ka P_2$ has nine words. But for every such $C_Q$ if $Q\cup \{v,w\}$ is a twin pair free code with four words, $v,w\in S^5$, and $C_Q\cap (Q\cup \{v,w\})=\emptyset$, then $C_Q$ is not a cover of  $Q\cup \{v,w\}$. 


For every $Q\in \ka P_1$ such that $|C_Q|=9$ it can be computed that if $Q\cup \{v,w\}$ is a twin pair free code having four words which is disjoint with $C_Q$, $v,w\in S^5$, then $C_Q$ is not a cover of  $Q\cup \{v,w\}$. For $Q\in \ka P_1$ if $C_Q$ contains eight words, then the  computations show that for every $v,w,u\in S^d$ if $C_Q\cup \{u\}$, $Q\cup \{v,w\}$ are disjoint twin pair free codes containing nine and four words, respectively, then $Q\cup \{v,w\}$ is not covered by $C_Q\cup \{u\}$. Similarly, if $|C_Q|=7$, then for every $v,w,u,q\in S^5$ if $C_Q\cup \{u,q\}$, $Q\cup \{v,w\}$ are disjoint twin pair free codes having nine and four words, respectively, then $Q\cup \{v,w\}$ is not covered by $C_Q\cup \{u,q\}$.   
\hfill{$\square$}

\vspace{-4mm}
\subsection{$q$-equivalent codes}

From Lemma \ref{po7}$i$, \ref{twin} and (\ref{dkostki}) we obtain 
\begin{wn}
\label{wa1}
For every $d\geq 1$ and every $q\in S^d$ there are no $q$-equivalent twin pair free disjoint codes $K,M\subset S^d$ such that $|K|=1$ and $|M|\in [4]$. 
\hfill{$\square$}
\end{wn}

\medskip
The following $q$-equivalent codes will be used throughout the paper: 

\begin{lemat}
\label{cody}
Let $K,M\subset S^d$ be  disjoint and twin pair free $q$-equivalent codes. 

\smallskip
$(i)$ If $|K|=1$ and $|M|=5$, then, up to isomorphism, $q_A=bbb$, $K_A=\{bbb\}$ and  $M_A=\{aaa,a'a'a',baa',a'ba,aa'b\}$,
where $A=\{1,2,3\}$, $K_{A^c}=M_{A^c}=\{b...b\}$ and $q_{A^c}=s_4...s_d$, $s_i\neq b'$ for $i\in \{4,...,d\}$. If $|K|=1$ and $|M|=6$, then, up to isomorphism, $q_A=bbbb$, $K_A=\{bbbb\}$ and  $M_A=\{aaaab,\;a'a'a'ab,\;baa'ab,\;a'baab,\;aa'bab,\;bbba'b\}$,
where $A=\{1,2,3,4\}$, $K_{A^c}=M_{A^c}=\{b...b\}$ and $q_{A^c}=s_5...s_d$, $s_i\neq b'$ for $i\in \{5,...,d\}$.


\smallskip
$(ii)$ If $|K|=|M|=2$, then, up to isomorphism, $q_A=bb$, $K_A=\{ab,a'a\}$ and  $M_A=\{ba,aa'\}$,
where $A=\{1,2\}$, $K_{A^c}=M_{A^c}=\{b...b\}$ and $q_{A^c}=s_3...s_d$, $s_i\neq b'$ for $i\in \{3,...,d\}$. 

\smallskip
$(iii)$ If $|K|=2,|M|=3$, then, up to isomorphism, $q_A=bbb$ and  
$$
K_A=\{abb, a'a'a'\},\;  M_A=\{ba'a', aba, aaa'\}\; or \; K_A=\{aab, ba'a'\},\;  M_A=\{aba', aaa, a'a'a'\},
$$
where $A=\{1,2,3\}$, $K_{A^c}=M_{A^c}=\{b...b\}$ and $q_{A^c}=s_4...s_d$, $s_i\neq b'$ for $i\in \{4,...,d\}$. 

\smallskip
$(iv)$ In $(i)$-$(iii)$ if $r_A\neq q_A$, where $r\in S^d$, then $K$ and $M$ are not $r$-equivalent.

\smallskip
$(v)$ There are, up to isomorphism, seventeen disjoint and twin pair free $q$-equivalent codes $K,M$ for $q=bbb$ written down in the alphabet $S=\{a,a',b\}$. These are

\vspace{-6mm}
\begin{center}
\begin{displaymath}
\begin{tabular}{|l|l|l|}
\hline
 $1$& $q=bbb$ & $K=\{aab,aa'a\},\;\; M=\{aba,aaa'\}$ \\ \cline{1-3}
 $2$& $q=bbb$ & $K=\{bab,ba'a\},\;\; M=\{bba,baa'\}$ \\ \cline{1-3}
 $3$& $q=bbb$ & $K=\{abb,\; a'a'a'\},\; M=\{ba'a',aba,aaa'\}$ \\ \cline{1-3}
 $4$& $q=bbb$ & $K=\{aab, ba'a'\},\; M=\{aba', aaa, a'a'a'\}$ \\ \cline{1-3}
 $5$& $q=bbb$ & $K=\{bab, ba'a\},\; M=\{aab, aa'a, a'ba, a'aa'\}$ \\ \cline{1-3}
 $6$& $q=bbb$ & $K=\{aba', a'a'a', bba\},\; M=\{aab, a'aa, ba'b\}$ \\ \cline{1-3}
 $7$& $q=bbb$ & $K=\{aba', a'a'a', bba\},\; M=\{abb, a'aa, a'a'b\}$ \\ \cline{1-3}
 $8$& $q=bbb$ & $K=\{aba', a'a'a', bba\},\; M=\{aaa', baa, ba'b\}$ \\ \cline{1-3}
 $9$& $q=bbb$ & $K=\{aab, a'ba', ba'a, \},\; M=\{aba, a'a'b, baa'\}$ \\ \cline{1-3}
 $10$& $q=bbb$ & $K=\{aba, a'aa', ba'a'\},\; M=\{aaa, aa'b, a'ba'\}$ \\ \cline{1-3}
 $11$& $q=bbb$ & $K=\{aaa', aba, a'aa, a'ba'\},\; M=\{aa'a, a'a'a', bab\}$ \\ \cline{1-3}
 $12$& $q=bbb$ & $K=\{aaa, aba', a'ab, ba'a\},\; M=\{aa'a', baa', bba\}$ \\ \cline{1-3}
 $13$& $q=bbb$ & $K=\{aba', a'aa', a'a'a, baa\},\; M=\{ aaa, aa'a', a'ba, baa'\}$ \\ \cline{1-3}
 $14$& $q=bbb$ & $K=\{aab, aa'a, a'aa', a'ba\},\; M=\{aaa', aba, a'ab, a'a'a\}$ \\ \cline{1-3}
 $15$& $q=bbb$ & $K=\{aaa, aba', a'a'a', a'ba\},\; M=\{aab, aa'a', a'aa, a'a'b\}$ \\ \cline{1-3}
 $16$& $q=bbb$ & $K=\{aaa, a'a'a', baa', a'ba,aa'b\},\;M=\{bbb\}$ \\ \cline{1-3}
 $17$& $q=bbb$ & $K=\{aaa, a'a'a', baa', a'ba,aa'b\},\; M=\{a'a'a, aaa', ba'a', aba,a'ab\}$ \\ \cline{1-3}
    \hline
\end{tabular}
\end{displaymath}
\end{center}

\vspace{-2mm}
\noindent{\footnotesize {\rm Table 2:  $q$-equivalent twin pairs-free disjoint codes $K,M$ for $q=bbb$. 
Recall that two $q$-equivalent codes $K,M$ are isomorphic to $q$-equivalent codes $K^1,M^1$ if there is an isomorphism $f$ such that $f(q)=q$ and $f(K)=K^1$ and $f(M)=M^1$.}
}
\end{lemat}

\vspace{-3mm}
\proof The part $(i)$ follows from the forms of the covers $C^1$ and $C^2$ given in Table 1 and the definition of $q$-equivalent codes.  

The simplest way to find all, up to isomorphism, $q$-equivalent twin pair free disjoint codes $K,M$ for $q=bbb$, which are enumerated in Table 2, is to compute the family $\ka C_3$ of all covers of the word $bbb$ having $n\in \{3,...,8\}$ words. 
Having $\ka C_3$ for every $C_1,C_2\in \ka C_3$ if $C_1\setminus C_1\cap C_2$ and $C_2\setminus C_1\cap C_2$ are disjoint and twin pair free codes, then $K=C_1\setminus C_1\cap C_2$ and $M=C_2\setminus C_1\cap C_2$. (Clearly, codes in Table 2 can be found by the trial and error method.)

Since it is easy to show that the codes given in $(ii)$ and $(iii)$ are such that the set $K_{A^c}$ is a singleton whose element belongs to $S^{d-2}$ for ($ii$) and to $S^{d-3}$ for ($iii$), it follows that from $\ka C_3$ we obtain also codes described in $(ii)$ and $(iii)$.

We prove the statement $(iv)$ for the case $(ii)$. (Along the same lines we prove this statement for $(i)$ and $(iii)$). Let us consider the following boxes: $K^1=(Eb'\cap Ea)\times (Ea\cap Eb)\times G$, $K^2=(Eb'\cap Ea')\times Ea\times G$, $K^3=(Eb\cap Ea)\times (Ea\cap Eb')\times G$ and $K^4=Ea\times (Ea'\cap Eb')\times G$, where $G=(Eb)^{d-2}$. Clearly, $K^1\cup K^2=\bigcup E(K)\setminus \bigcup E(M)$ and $K^3\cup K^4=\bigcup E(M)\setminus \bigcup E(K)$. It follows from the definition of $q$-equivalent codes that for every $i\in [4]$ the box $K^i$ cannot be intersected by the box $\breve{q}$. Clearly, if $q_i\in \{a,a',b'\}$ for some $i\in \{1,2\}$, then $\breve{q}$ does not intersect all boxes in $K$ or in $M$, which is impossible, by Corollary \ref{wa1}. Now, by (\ref{dkostki}), it is easy to check that if $q$ is such that $q_i\in S\setminus \{a,a',b,b'\}$ for some $i\in \{1,2\}$, then $\breve{q}\cap K^j\neq\emptyset$ for some $j\in [4]$, which is impossible. 
\hfill{$\square$}

\section{Distribution of words in equivalent codes}
In general, for equivalent codes $V$ and $W$ we have $D_i(V)\neq D_i(W)$, $i\in [d]$, (see Subsection 2.3 for the definition of $D_i(V)$). However, in some cases to find the structure of equivalent codes $V$ and $W$ it is useful to compare the pairs  $(V^{i,l},V^{i,l'})$  and $(W^{i,l},W^{i,l'})$, $i\in [d], l,l'\in S$. In this section we describe selected relationships between  $(V^{i,l},V^{i,l'})$  and $(W^{i,l},W^{i,l'})$.

\smallskip
In \cite{Kis22} we proved the following two lemmas (Lemma 5 and Lemma 13 of \cite{Kis22}):
\begin{lemat}
\label{cy}
Let $S=\{a_1,a_1',...,a_k,a_k'\}$, $\varepsilon \in \{0,1\}^k$ and $A_\varepsilon=\{a_1^{\varepsilon_1},...,a_k^{\varepsilon_k}\}$, where $a_n^0=a_n,a^1_n=a'_n$ for $n\in [k]$. If $V,W\subseteq S^d$ are equivalent polybox codes, then for every $i\in [d]$ and every $\varepsilon \in \{0,1\}^k$ the codes $\bigcup_{l\in A_\varepsilon}V_{i^c}^{i,l}$ and $\bigcup_{l\in A_\varepsilon}W_{i^c}^{i,l}$
are equivalent.
\hfill{$\square$}
\end{lemat}

\begin{lemat}
\label{po1m}
Let $V,W\subset S^d$ be equivalent codes. Suppose that $w\in W^{i,s}\cup W^{i,s'}$ is a word such that the sets
$$
U^{i,l}=\{v \in V^{i,l}\colon \breve{w}\cap \breve{v}\neq\emptyset\},\;\; U^{i,l'}=\{v \in V^{i,l'}\colon \breve{w}\cap \breve{v}\neq\emptyset\}
$$
are non-empty, where $l\not\in \{s,s'\}$. If $|U^{i,l}|=1$, then $U^{i,l}_{i^c}\sqsubseteq U^{i,l'}_{i^c}$.
\hfill{$\square$}
\end{lemat}

\medskip
From  Corollary \ref{wa1} and Lemma \ref{wac} we obtain
	
\begin{lemat}
\label{wa}
If $V,W\subset S^d$ are disjoint equivalent polybox codes, the code $V$ does not contain a twin pair and there are $i\in [d]$ and $l\in S$ such that $|V^{i,l}|=1$ and $|V^{i,l'}|\in [4]$, then $V^{i,l}\sqsubseteq W^{i,l}$ and  $V^{i,l'}\sqsubseteq W^{i,l'}$. 
\hfill{$\square$}
\end{lemat}

\begin{lemat}
\label{r24}
Let $V,W\subset S^d$ be disjoint equivalent polybox codes without twin pairs. If there are $i\in [d]$ and $l\in S$ such that $V=V^{i,l}\cup V^{i,l'}$  and $V^{i,l}\neq\emptyset$, $V^{i,l'}\neq\emptyset$, then $|V|\geq 24$.
\end{lemat}
\proof If $W=W^{i,l}\cup W^{i,l'}$, then the codes $V^{i,l},W^{i,l}$ are disjoint and equivalent and similarly, $V^{i,l'},W^{i,l'}$ are disjoint and equivalent.  Thus, by Theorem \ref{t12}, $|V^{i,l}|\geq 12$ and $|V^{i,l'}|\geq 12$. 
 
If $W^{i,s}\neq\emptyset$ for some $s\not\in \{l,l'\}$, then, by ({\bf C}), $W^{i,s}_{i^c}$ and  $W^{i,s'}_{i^c}$ are equivalent. Since these codes are disjoint and twin pair free, again by Theorem \ref{t12}, $|W_{i^c}^{i,s}|\geq 12$ and $|W_{i^c}^{i,s'}|\geq 12$. 
\hfill{$\square$}

\medskip
Recall that a code $V\subset S^d$ is {\it flat} if there are $i\in [d]$ and $l\in S$ such that $V=V^{i,l}$.

\begin{lemat}
\label{rr2}
Let $V,W\subset S^d$, where $S=\{a,a',b,b',c,c'\}$ and $d=5,6$, be disjoint equivalent polybox codes without twin pairs, $|V|\in \{13,...,16\}$, and let  
$V$ be not flat. For every $i\in [d]$ and every $l\in S$, if $V^{i,l}\cup V^{i,l'}\neq\emptyset$, then  $V^{i,l}\neq\emptyset$ and $V^{i,l'}\neq\emptyset$. 
\end{lemat}
\proof Suppose on the contrary that there are $i\in [d]$ and $l\in S$ such that $V^{i,l'}=\emptyset$ and $V^{i,l}\neq\emptyset$. Then, by ({\bf P}), $V^{i,l}\sqsubseteq W^{i,l}$. We may assume that $l=a$.


If the codes $V^{i,a}$ and $W^{i,a}$ are equivalent, then, since they are disjoint and twin pair free, by Theorem \ref{t12}, $|V^{i,a}|\geq 12$. Observe now that the equivalence of the codes $V^{i,a}, W^{i,a}$ implies the equivalence of the disjoint and twin pair free codes $V\setminus V^{i,a}$ and $W\setminus W^{i,a}$. Then, again by Theorem \ref{t12}, $|V\setminus V^{i,a}|\geq 12$, and thus $|V|\geq 24$, a contradiction. 

It follows from the above that we can assume that $V^{i,l}$ and $W^{i,l}$ are not  equivalent for $l\in S$.

Since the codes $V^{i,a}$ and $W^{i,a}$ are not equivalent, we have $W^{i,a'}\neq\emptyset$. As $V^{i,a'}=\emptyset$, it follows that, by Lemma \ref{rr1}, $W_{i^c}^{i,a'}\sqsubseteq W_{i^c}^{i,a}$ and the set $W_{i^c}^{i,a'}\cup V_{i^c}^{i,a}$ is a polybox code. In other words, taking into account that  $V^{i,a}\sqsubseteq W^{i,a}$, the codes $W^{i,a}_{i^c}$ and $W_{i^c}^{i,a'}\cup V_{i^c}^{i,a}$ are equivalent. Thus, $|W^{i,a}|-|W^{i,a'}|=|V^{i,a}|$ and, by Lemma \ref{po7}$iii$, $|W^{i,a}|\geq 10$. Hence,   
\begin{equation}
\label{k16}
|W^{i,a}|+|W^{i,a'}|\geq 11. 
\end{equation}

Assume first that $V^{i,b}\cup V^{i,b'}\neq\emptyset$ and $V^{i,c}\cup V^{i,c'}=\emptyset$. 
It follows that, $W^{i,c}\cup W^{i,c'}=\emptyset$, otherwise,  by ({\bf C}), the codes $W^{i,c}_{i^c}, W^{i,c'}_{i^c}$ are equivalent, and thus, by Theorem \ref{t12}, $|W^{i,c}_{i^c}|, |W^{i,c'}_{i^c}|\geq 12$, a contradiction. 
Then $V^{i,l}\neq\emptyset$ for $l\in \{b,b'\}$, as if it is not so, in the same manner as (\ref{k16}), we show that $|W^{i,b}|+|W^{i,b'}|\geq 11$, and consequently $|W|\geq 22$, which is not true. Since $|W^{i,a}|\geq 10$, we have $|W^{i,a'}|+|W^{i,l}|\leq 6$ for $l\in \{b,b'\}$. 
By Lemma \ref{cy}, $|W^{i,a'}|+|W^{i,l}|=|V^{i,a'}|+|V^{i,l}|$ for $l\in \{b,b'\}$, and thus $|V^{i,l}|\leq 6$ for $l\in \{b,b'\}$. 

{\sc Claim} The code $V^{i,l}$ is rigid.

To prove the claim let $V^{i,l}=Q$. If $|Q|\leq 5$, then the rigidity of $Q$ follows from Lemma \ref{po7}$i$ and Lemma 9 of \cite{Kis22}. Let $|Q|=6$ and assume on the contrary that there is a code $W\subset S^d$ which is equivalent to $Q$ and disjoint with it.
By Theorem \ref{t12}, there is a twin pair in $W$. Let $u,w\in W$ be such a pair. We may assume that $u_1=w_1'$ and $u_{1^c}=w_{1^c}$. Let $K=\breve{u}\cup \breve{w}$. Since $K\subset \bigcup E(Q)$, the family of boxes $\ka F=\{K\cap \breve{v}\colon v\in Q\}$ is a suit for $K$ which, by Lemma \ref{twin}, does not contain a twin pair. Let $\ka F^{1,s}=\{K\cap \breve{v}\colon v\in Q^{1,s}\}$, $s\in S$. Note that $|\ka F^{1,s}|=5$ and $|\ka F^{1,s'}|=1$, that is $|Q^{1,s}|=5$ and $|Q^{1,s'}|=1$. By Lemma \ref{po1m}, $Q^{1,s'}_{1^c}\sqsubseteq Q^{1,s}_{1^c}$, and thus $Q$ may be written down in the alphabet $\{a,a',b,b'\}$ (see the code $C^1$ in Table 1). In Lemma 16 of \cite{Kis22} we showed that every code $Q\subset \{a,a',b,b'\}^d$ with $|Q|\leq 9$ is rigid. The claim has been proved

By Lemma \ref{cy}, the codes $V^{i,l}_{i^c}$ and $W^{i,a'}_{i^c}\cup W^{i,l}_{i^c}$ are equivalent for $l\in \{b,b'\}$. As $V\cap W=\emptyset$ and the code $V^{i,l}_{i^c}$ is rigid for $l\in \{b,b'\}$, we have  $V^{i,l}_{i^c}=W^{i,a'}_{i^c}$ for $l\in \{b,b'\}$. Then  $V^{i,b}\cup V^{i,b'}$ contains a twin pair. A contradiction.

Therefore, we may assume that $V^{i,l}\neq\emptyset$ for $l\in \{b,b',c,c'\}$. Then, as before, we show that $W^{i,l}\cup W^{i,l'}\neq\emptyset$ for $l\in \{b,c\}$. 
Note that we may assume that $W^{i,l}\neq\emptyset$ for $l\in \{b,b'\}$, because if on the contrary $W^{i,b}=W^{i,c}=\emptyset$, then, in the same manner as (\ref{k16}), we show that $|V^{i,b}|+|V^{i,b'}|\geq 11$ and  $|V^{i,c}|+|V^{i,c'}|\geq 11$ which is not possible.


Observe that  $|V^{i,a}|\geq 5$, otherwise $|W^{i,a'}|\geq 6$, and then $W=W^{i,a}\cup W^{i,a'}$ which, by Lemma \ref{r24}, means that $|W|\geq 24$. A contradiction. If $W^{i,c}=\emptyset$, then $|V^{i,c}|+|V^{i,c'}|\geq 11$, and consequently $V^{i,b}\cup V^{i,b'}=\emptyset$ which is not true.
Therefore, $W^{i,l}\neq\emptyset$ for every $l\in \{b,b',c,c'\}$. By (\ref{k16}), we may assume that $|W^{i,b}|+|W^{i,b'}|\leq 3$. By Lemma \ref{wa}, $W^{i,l}\sqsubseteq V^{i,l}$ for $l\in \{b,b'\}$, and by Lemma \ref{po7}$i$,  $|V^{i,b}|+|V^{i,b'}|\geq 10$. But then $V^{i,c}=\emptyset$ or $V^{i,c'}=\emptyset$, a contradiction.
\hfill{$\square$}

\begin{lemat}
\label{rr3}
Let $V,W\subset S^d$, where $S=\{a,a',b,b',c,c'\}$ and $d=5,6$, be two disjoint equivalent polybox codes without twin pairs, $|V|\in \{13,...,16\}$, and let $V$ be not flat. Then

\smallskip
$(i)$ $\{|V^{i,l}|,|V^{i,l'}|\}\neq \{1,M\}$, where $M\in \{2,3,4\}$, and $|V^{i,l}|\leq 11$ for every $i\in [d]$ and $l\in S$.

\smallskip
$(ii)$ If $|V^{i,l}|\geq 2$ and $|V^{i,l'}|\geq 2$  or $|V^{i,l}|=1$ and $|V^{i,l'}|\geq 5$ 
for some $i\in [d]$ and $l\in S$, then there is $w\in W^{i,s}\cup W^{i,s'}$, $s\not \in \{l,l'\}$, such that $\breve{w}\cap \bigcup E(V^{i,l}\cup V^{i,l'})\neq\emptyset$.
\end{lemat}

\proof Fix $i\in [d]$. We first prove the lemma in the case $V^{i,l}\cup V^{i,l'}\neq\emptyset$ for every $l\in \{a,b,c\}\subset S$, which is an easy task.  Observe that $W^{i,l}\cup W^{i,l'}\neq\emptyset$ for $l\in \{a,b,c\}$, otherwise, by ({\bf C}), the codes $V^{i,l}_{i^c}, V^{i,l'}_{i^c}$ are equivalent for some $l\in \{a,b,c\}$, and thus, by Theorem \ref{t12}, $|V^{i,l}\cup V^{i,l'}|\geq 24$, a contradiction. By Lemma \ref{rr2}, $V^{i,l},W^{i,l}\neq\emptyset$ for $l\in \{a,a',b,b',c,c'\}$.

We may assume on the contrary that $|V^{i,a}|=1$ and $|V^{i,a'}|=M$ for some $M\in \{2,3,4\}$. Then, by Lemma \ref{wa}, \ref{po7}$i$ and \ref{po7}$ii$, $|W^{i,a}\cup W^{i,a'}|= 12$ if $M=2$ and $|V|=16$, and thus $|W^{i,l}|=1$ for $l\in \{b,b',c,c'\}$; if $M=2$ and $|V|\leq 15$ or $M\in \{3,4\}$ and $|V|\leq 16$, then $W^{i,l}=\emptyset$ for some $l\in \{b,b',c,c'\}$ . In the first case, again by Lemma \ref{wa} and \ref{po7}$i$, $|V^{i,l}|\geq 5$ for $l\in \{b,b',c,c'\}$ which contradicts the assumption $|V|\leq 16$. The second case is not possible, by Lemma \ref{rr2}.

If $|V^{i,l}|\geq 12$, then $V^{i,s}=\emptyset$ for some $s\in \{a,a',b,b',c,c'\}\setminus \{l\}$, which is a contradiction, by Lemma \ref{rr2}.

Finally, if $|V^{i,l}|\geq 2$ and $|V^{i,l'}|\geq 2$ or $|V^{i,l}|=1$ and $|V^{i,l'}|\geq 5$ and in both cases 
 on the contrary  $V^{i,l}\sqsubseteq W^{i,l}$, $V^{i,l'}\sqsubseteq W^{i,l'}$ , then, by Lemma \ref{po7} and \ref{rr2}, $W^{i,s}\cup W^{i,s'}=\emptyset$ for some  $s\in \{a,b,c\}\setminus \{l,l'\}$, which is, as we showed, impossible. This completes the proof of the lemma for the case $V^{i,l}\cup V^{i,l'}\neq\emptyset$ for $l\in \{a,b,c\}\subset S$.

\smallskip
Therefore, in what follows we may assume that $V^{i,l}\cup V^{i,l'}=\emptyset$ for  $l\in \{c,c'\}$

\smallskip 
To prove the second part of $(i)$ assume on the contrary that $|V^{i,l}|\geq 12$ for some $l\in S$. Then, by Lemma \ref{r24} and \ref{rr2}, $|V^{i,s}|=|V^{i,s'}|=1$ or $|V^{i,s}|=1$, $|V^{i,s'}|=2$, where $s\not\in \{l,l'\}$. Therefore, by Lemma \ref{wa} and \ref{po7}$i$, $|W^{i,s}\cup W^{i,s'}|\geq 10$ from where  $|W^{i,l}\cup W^{i,l'}|\leq 6$. Consequently, $m(W^{i,l}_{i^c}\setminus W^{i,l'}_{i^c})\leq 6$. On the other hand $m(V^{i,l}_{i^c}\setminus V^{i,l'}_{i^c})\geq 10$ which contradicts (\ref{mb}).  

\smallskip
To prove the first part of $(i)$ suppose on the contrary that  $|V^{i,l}|=1$ and $|V^{i,l'}|=M$ for some $l\in S$, where $M\in \{2,3,4\}$. Then, by Lemma \ref{wa}, $V^{i,l}\sqsubseteq W^{i,l}$ and $V^{i,l'}\sqsubseteq W^{i,l'}$.

Let $|V|=16$, $M=2$ and $d=6$.  By Lemma \ref{po7}$i$ and \ref{po7}$ii$, $|W^{i,l}|\geq 5$ and $|W^{i,l'}|\geq 7$. Taking into account Lemma \ref{cy}, \ref{po7} and (\ref{mb}), one can easily show that  $|W^{i,l}|=6$,$|W^{i,l'}|=7$, $|W^{i,s}|=1$,$|W^{i,s'}|=2$ and $|V^{i,s}|=6$,$|V^{i,s'}|=7$. The relations $V^{i,l}\sqsubseteq W^{i,l}$ and $V^{i,l'}\sqsubseteq W^{i,l'}$ allow us to partially predict the structures of the codes  $W^{i,l},W^{i,l'}$. 

Namely, if $W^{i,l}$ contains six words and $\breve{v}\cap \breve{w}\neq\emptyset$ for $v\in V^{i,l}$ and for every $w\in W^{i,l}$, then, up to isomorphism,  $W^{i,l}=\{wb\colon w\in C^2\}$, where $C^2$ is given in Table 1; if  $W^{i,l}$ contains six words but $\breve{v}\cap \breve{w}\neq\emptyset$ for five words $w\in W^{i,l}$, then, up to isomorphism,  $W^{i,l}=\{wb\colon w\in C^1\}\cup \{u\}$, where $C^1$ is given in Table 1 and $u\in S^d$. Finally, by Lemma \ref{po7}$ii$, we may assume that $V^{i,l'}=\{l'bbbbb,\;l'b'abbb\}$ and $W^{i,l'}_{i^c}$ is equal to the code $C_2$ given in that lemma. 

The computations show  that there is only one, up to isomorphism, pair of the codes  $V^{i,l}, V^{i,l'}$ and $W^{i,l},W^{i,l'}$ with the above properties for which (\ref{dif}) is valid (we took $i=1$ and $l=b$):
$$
V^{1,b}=\{ba'aab'a\},\;\; V^{1,b'}=\{b'bbbbb,\;b'b'abbb\} 
$$
and
$$
W^{1,b}=\{bbb'bbb,ba'aaab', ba'aba'b', ba'abb'b, ba'ab'ab, ba'ab'a'a\},
$$
$$
W^{1,b'}=\{b'ba'bbb, b'aabbb, b'a'aaaa', b'a'aaba, b'a'aa'ab, b'a'aa'a'a, b'a'aba'a'\}. 
$$
Since $D_3(W^{1,b}\cup W^{1,b'})=((11,1),(0,1), (0,0))$, in the similar way as at the beginning of the proof of ($i$), we show that $D_3(W)=((11,1),(2,2),(0,0))$. 
Clearly, there is $v\in V\setminus (V^{3,a}\cup V^{3,a'})$ such that $\breve{v}\cap \bigcup E(W^{3,a}\cup W^{3,a'})\neq\emptyset$, otherwise  $W^{3,l}\sqsubseteq V^{3,l}$ for $l\in \{a,a'\}$, and  consequently, by Theorem \ref{t12} and Lemma \ref{po7}$i$, $|V^{3,a}|\geq 12$ and $|V^{3,a'}|\geq 5$, which is not possible. 

Let $U=U^{3,a}\cup U^{3,a'}$, where $U^{3,a}=W^{1,b'}\setminus \{b'ba'bbb\}$ and $U^{3,a'}=\{b'ba'bbb\}$. We shall use Lemma \ref{ci} in which $i=3$ and $l=a$ but we change the role of $V$ and $W$. We have $U^{3,l}\subseteq W^{3,l}$ for $l\in \{a,a'\}$. Observe that, $\breve{w}_{3^c}\cap \breve{u}_{3^c}=\emptyset$ for every $w\in W^{3,a}\setminus U^{3,a}$, where $\{u\}=U^{3,a'}$. For this reason, if $v\in V\setminus (V^{3,a}\cup V^{3,a'})$ is such that $\breve{v}\cap \breve{w}\neq\emptyset$, where $w\in W^{3,a}\cup W^{3,a'}$, then $w\in U$. Recall that for every $v\in V\setminus (V^{3,a}\cup V^{3,a'})$ the set $\breve{v}\cap \bigcup E(U)$ is a $3$-cylinder in $\breve{v}$. Therefore, for every $v\in V\setminus (V^{3,a}\cup V^{3,a'})$ such that $\breve{v}\cap \bigcup E(U)\neq\emptyset$ we have $v_A=u_A=bbbb$, where $A=\{2,4,5,6\}$. Moreover, $U_B=\{b'\}$, where $B=\{1\}$. By Lemma \ref{ci}, in which $p=bbbb$ and $q=b'$, the code $V$ or $W$ contains a twin pair. A contradiction. 
This completes the proof of the case $\{|V^{i,l}|,|V^{i,l'}|\}\neq \{1,2\}$ for $d=6$ and $|V|=16$

\smallskip
If $d=5$, $M=2$ and $|V|$=16, then, by Lemma \ref{po7}$i$ and \ref{po7}$iii$,  $|W^{i,l}|\geq 5$ and $|W^{i,l'}|\geq 8$. By (\ref{mb}), 
must be $|W^{i,l}|\geq 7$. Indeed, if $|W^{i,l}|\leq 6$, then, to preserve the equality (\ref{mb}), $m(W^{i,l'}_{i^c}\setminus W_{i^c}^{i,l})=2$ and consequently, $W^{i,l}_{i^c}\sqsubseteq  W_{i^c}^{i,l'}$. Then, by Lemma \ref{po7}$iii$, $|W^{i,l'}|\geq 10$. 
Therefore, by Lemma \ref{rr2}, $W=W^{i,l}\cup W^{i,l'}$ and then $|W|\geq 24$, by Lemma \ref{r24}. A contradiction.  Thus, $|W^{i,l}|\geq 7$, which means that $|W|>16$, because $W^{i,s},W^{i,s'}\neq\emptyset$ for $s\in \{a,a',b,b'\}\setminus \{l,l'\}$. A contradiction.
 
\smallskip
Let now $M\in \{3,4\}$, $d=5,6$ and $|V|=16$. By Lemma \ref{po7}, $|W^{i,l}|\geq 5$ and $|W^{i,l'}|\geq 9$. 
By Lemma \ref{r24} and \ref{rr2}, we have $|W^{i,l}|=5$ and $|W^{i,l'}|=9$. Since $m(W^{i,l'}_{i^c}\setminus W^{i,l}_{i^c})\geq 4$ and $m(V^{i,l'}_{i^c}\setminus V^{i,l}_{i^c})\leq 4$, by (\ref{mb}), $m(W^{i,l'}_{i^c}\setminus W^{i,l}_{i^c})=4$ which means that $m(W^{i,l}_{i^c}\setminus W^{i,l'}_{i^c})=0$. Thus, $V^{i,l}_{i^c}\sqsubseteq V^{i,l'}_{i^c}$. A contradiction, by Lemma \ref{po7}$i$. The proof of ($i$) is completed for $|V|=16$.



Let $|V|=15$ and $d=6$.

If $M=2$, then, by (\ref{mb}), Lemma \ref{po7}$ii$, \ref{rr2} and \ref{cy}, it is easy to see that $|W^{i,l}|=6$ and $|W^{i,l'}|=7$, $|W^{i,s}|=|W^{i,s'}|=1$ and $|V^{i,s}|=|V^{i,s'}|=6$. In the similar way as before (the case $M=2, |V|=16, d=6$) we compute all pairs of codes $W^{i,s},W^{i,s'}$ and $V^{i,s},V^{i,s'}$ for which the last two equalities and (\ref{dif}) are valid. The computations show that there are no codes $W^{i,s},W^{i,s'}$ and $V^{i,s},V^{i,s'}$ with the above properties, where additionally, by the first part of the proof, $|(V^{i,s}\cup V^{i,s'})^{j,l}|\leq 11$ for every $j\in [d]\setminus \{i\}$ and $l\in S$. 

If $M=3$, then $|W^{i,l}|\geq 5$, $|W^{i,l'}|\geq 9$, and if $M=4$, then $|W^{i,l}|\geq 5$, $|W^{i,l'}|\geq 10$, by Lemma \ref{po7}. In both cases the set $W^{i,s}$ or $W^{i,s'}$ is empty, which is a contradiction, by Lemma \ref{rr2} and \ref{r24}. Along the same lines we prove $(i)$ for  $|V|=15$ and $d=5$ as well as for $|V|\in \{13,14\}$ and $d=5,6$

\smallskip
To prove $(ii)$, let $|V^{i,l}|=|V^{i,l'}|=2$ and  suppose on the contrary that $\bigcup E(W^{i,s}\cup W^{i,s'})\cap \bigcup E(V^{i,l}\cup V^{i,l'})=\emptyset$ for every $s\in S\setminus \{l,l'\}$. Then $V^{i,l}\sqsubseteq W^{i,l}$ and $V^{i,l'}\sqsubseteq W^{i,l'}$.

Assume first that $d=6$. By Lemma \ref{po7}$ii$, $|W^{i,l}|\geq 7$ and $|W^{i,l'}|\geq 7$. 
By Lemma \ref{r24} and \ref{rr2}, $|W^{i,s}|=|W^{i,s'}|=1$ for some $s\in S\setminus \{l,l'\}$.  
Then, by Lemma \ref{rr1},  $|V^{i,s}|=|V^{i,s'}|=6$. This, as it was showed in the proof of the part ($i$), is not possible.

Let now $d=5$. Since $V^{i,l}\sqsubseteq W^{i,l}$ and $V^{i,l'}\sqsubseteq W^{i,l'}$, by Lemma \ref{po7}$iii$, we have $|W^{i,l}|\geq 8$ and $|W^{i,l'}|\geq 8$ and  consequently $W=W^{i,l}\cup W^{i,l'}$. Then, by Lemma \ref{r24}, $|W|\geq 24$, a contradiction.

For the same reason, if $|V^{i,l}|\geq 3$, $|V^{i,l'}|\geq 2$ and  $V^{i,l}\sqsubseteq W^{i,l}$, $V^{i,l'}\sqsubseteq W^{i,l'}$, then  $|W|\geq 24$. 

Finally, let $|V^{i,l}|=1$ and $|V^{i,l'}|\geq 5$ and $V^{i,l}\sqsubseteq W^{i,l}$, $V^{i,l'}\sqsubseteq W^{i,l'}$. Then, by Lemma \ref{po7}$iii$, $|W^{i,l}|\geq 5$ and $|W^{i,l'}|\geq 10$, and then, by Lemma \ref{rr2}, $W=W^{i,l}\cup W^{i,l'}$ which gives $|W|\geq 24$, a contradiction.
This completes the proof.
\hfill{$\square$}

\smallskip
At the end of this section we summarize the above results in the following lemma which is a base for reductions that will be made in the final computations. Recall that if $S=\{a_1,a_1',...,a_k,a'_k\}$, $Q\subset S^d$ and $i\in [d]$, then $D_i(Q)=((|Q^{i,a_1}|,|Q^{i,a_1'}|),...,(|Q^{i,a_k}|,|Q^{i,a_k'}|))$. 

\begin{lemat}
\label{rr5}
Let $V,W\subset S^d$, where $d=5,6$, be two disjoint equivalent polybox codes without twin pairs, $|V|\in \{13,...,16\}$, and let $V$ be not flat. If $U\subseteq V$, $C_U\subseteq W$ is a cover of $U$ and $Q\subseteq C_U$, then:

\medskip
$(i)$ If $S=\{a,a',b,b',c,c'\}$, then $n,m \leq 11$ for every $i\in [d]$ and every pair $(n,m)$ in $D_i(Q)$.  

\medskip
$(ii)$ If $|Q|=15$ and $S=\{a,a',b,b'\}$, then for every $i\in [d]$ and every pair $(n,m)$ in $D_i(Q)$ we have $n+m>0$ and $n\geq 1$ if $2\leq m\leq 4$ (or conversely with respect to $n$ and $m$).

\medskip
$(iii)$ If $|Q|=15$ and $S=\{a,a',b,b',c,c'\}$, then for every $i\in [d]$ and every two pairs $(n,m), (p,q)$ in $D_i(Q)$ if $n+m=1$, then $p, q\leq 7$.

\medskip
$(vi)$ If $|Q|=16$ and $S=\{a,a',b,b'\}$, then for every $i\in [d]$ and for every pair $(n,m)$ in $D_i(Q)$ we have $n\cdot m>0$ and $n\geq 2$ if $2\leq m\leq 4$ (or conversely with respect to $n$ and $m$).

\medskip
$(v)$ If $|Q|=16$ and $S=\{a,a',b,b',c,c'\}$, then for every $i\in [d]$ and every two pairs $(n,m), (p,q)$ in $D_i(Q)$ if $n=m=1$, then $p, q\leq 7$.
\end{lemat}
\proof
By Lemma \ref{r24}, \ref{rr2} and \ref{rr3}, only the third and the fifth  statement needs an explanation. By Lemma \ref{rr2}, we assume that $|V|=16$. Suppose on the contrary that there are $(n,m)$ and $(p,q)$ in $D_i(Q)$ for some $i\in [d]$ such that $n=1$, $m=0$ and $p\geq 8$. Then, since $|Q|=15$, adding a new word $v$ to $Q$ such that $Q\cup \{v\}$ is a twin pair free code, we obtain the pair $(n',m')$, belonging to $D_i(Q\cup \{v\})$, which arises from $(n,m)$ and $(n',m')=(n,m)$ or $n'=2$ and $m'=0$ or  $n'=m'=1$. The first two cases can be excluded by Lemma \ref{rr2}. In the third case, the distribution $D_i(Q\cup \{v\})$ contains the pair $((1,1),(p,q))$, where $p\geq 8$, and thus $q\leq 6$. Since $Q\cup \{v\}=W$, we may assume that $|W^{i,a}|=|W^{i,a'}|=1$ and $|W^{i,b}|\geq 8$ and $|W^{i,b'}|\leq 6$. By Lemma \ref{wa} and \ref{po7}$i$, $V^{i,l}\geq 5$ for $l\in \{a,a'\}$. 
Suppose that $|V^{i,b}|\geq 2$ and $|V^{i,b'}|\geq 2$ or  $|V^{i,b}|\geq 5$ and $|V^{i,b'}|=1$. Then, by Lemma \ref{rr3}$ii$, there are $w\in W^{i,a}\cup W^{i,a'}$ and $v\in V^{i,b}\cup V^{i,b'}$  such that $\breve{w}\cap \breve{v}\neq\emptyset$. Since the set $\{\breve{v}\cap \breve{w}\colon w\in W^{i,a}\cup W^{i,a'}\}$ is an $i$-cylinder in $\breve{v}$, it follows that $W_{i^c}^{i,a}$ and $W_{i^c}^{i,a'}$ are $v_{i^c}$-equivalent.
This means that the words in $W^{i,a}\cup W^{i,a'}$ form a twin pair. A contradiction. Thus, by Lemma \ref{rr3}$i$, $|V^{i,b}|=|V^{i,b'}|=1$. Then $m(V_{i^c}^{i,b}\setminus V_{i^c}^{i,b'})\leq 1$. On the other hand  $m(W_{i^c}^{i,b}\setminus W_{i^c}^{i,b'})\geq 2$, which is a contradiction, by (\ref{mb}). The fifth statement is showed in the same manner. 
\hfill{$\square$}

\section{A rigidity result}

In this section we first prove a rigidity result, and next we use it to eliminate certain configurations of words in $V$ and $W$ (Corollary \ref{le12}). 

\begin{lemat}
\label{rig}
Every polybox code $V\subset S^d$ without twin pairs, where $S=\{a,a',b,b'\}$ and $d=4,5$, having at most 10 words is rigid. 
\end{lemat}
\proof
It follows from \cite[Lemma 9]{Kis22} that the lemma is true for $d\leq 3$. We give the proof for $d=4$. (The proof for $d=5$ is nearly the same, however during the proof we give the differences between the case $d=4$ and $d=5$).

Assume on the contrary that there is a code $W\subset S^d$ which is equivalent to $V$ and disjoint with it. We now proceed as in the proof of Lemma \ref{rr2}: By Theorem \ref{t12}, there is a twin pair, say $u,w$, in $W$. We may assume that $u_1=w_1'$ and $u_{1^c}=w_{1^c}$. Moreover, let $K=\breve{u}\cup \breve{w}$ and $\ka F=\{K\cap \breve{v}\colon v\in V\}$. Recall that, by Lemma \ref{twin}, the suit $\ka F$ for $K$ does not contain a twin pair. Below we enumerate hypothetical structures of  $\ka F$. Let $\ka F^{1,l}=\{K\cap \breve{v}\colon v\in V^{1,l}\}$.

\begin{enumerate}
\item $|\ka F|=m,\;\;m\in \{6,...,10\},\;\; {\rm and}\;\; |\ka F^{1,a}|=m-1,\;\; |\ka F^{1,a'}|=1$
\item $|\ka F|=8\;\; {\rm and}\;\; |\ka F^{1,a}|=2,\;\; |\ka F^{1,a'}|=2,\;\; |\ka F^{1,b}|=2,\;\; |\ka F^{1,b'}|=2$
\item $|\ka F|=9\;\; {\rm and}\;\; |\ka F^{1,a}|=2,\;\; |\ka F^{1,a'}|=2,\;\; |\ka F^{1,b}|=2,\;\; |\ka F^{1,b'}|=3$
\item $|\ka F|=10\;\; {\rm and}\;\; |\ka F^{1,a}|=2,\;\; |\ka F^{1,a'}|=3,\;\; |\ka F^{1,b}|=2,\;\; |\ka F^{1,b'}|=3$
\item $|\ka F|=10\;\; {\rm and}\;\; |\ka F^{1,a}|=2,\;\; |\ka F^{1,a'}|=2,\;\; |\ka F^{1,b}|=2,\;\; |\ka F^{1,b'}|=4$
\item $|\ka F|=10\;\; {\rm and}\;\; |\ka F^{1,a}|=2,\;\; |\ka F^{1,a'}|=2,\;\; |\ka F^{1,b}|=3,\;\; |\ka F^{1,b'}|=3$
\item $|\ka F|=10\;\; {\rm and}\;\; |\ka F^{1,a}|=2,\;\; |\ka F^{1,a'}|=2,\;\; |\ka F^{1,b}|=1,\;\; |\ka F^{1,b'}|=5$
\item $|\ka F|=10\;\; {\rm and}\;\; |\ka F^{1,a}|=5,\;\; |\ka F^{1,a'}|=5.$
\end{enumerate}

To decide whether the above cardinalities of $\ka F^{1,l}$, $l\in S$, are possible we first show that for every $i\in [d]$ and every $l\in S$ the sets $V^{i,l},W^{i,l}\subset S^d$ are nonempty. Suppose on the contrary that there is a letter $l\in S$, say $l=b$, such that $V^{i,b}=\emptyset$ for some $i\in [d]$. We may assume $V^{i,a}\neq\emptyset$. By Lemma \ref{cy}, the codes $V^{i,a}_{i^c}$ and $W^{i,a}_{i^c}\cup W^{i,b}_{i^c}$ are equivalent. The code $V^{i,a}_{i^c}\subset S^{3}$ is rigid and therefore $V^{i,a}_{i^c}=W^{i,b}_{i^c}$, for otherwise $V^{i,a}\cap W^{i,a}\neq\emptyset$. Then, by ({\bf Co}), $V^{i,a}_{i^c}\sqsubseteq V^{i,a'}_{i^c}$. In the same way we show that $V^{i,a'}_{i^c}\sqsubseteq V^{i,a}_{i^c}$, and consequently $V^{i,a}_{i^c}$ and $V^{i,a'}_{i^c}$ are equivalent which gives $|V^{i,a}_{i^c}|\geq 12$, by Theorem \ref{t12}. A contradiction. Therefore, $V^{i,l}\neq\emptyset$ for $l\in S$.

Now assume that $W^{i,b}=\emptyset$ and $|W^{i,b'}|=k>0$. Then, by ({\bf P}), $W^{i,b'}\sqsubseteq  V^{i,b'}$.

If $k\geq 3$, then $W^{i,b'}$ contains at least two words which do not form a twin pair, and thus, by Lemma \ref{po7}$ii$, $|V^{i,b'}|\geq 7$. Then, $|V^{i,b'}|=7$, $|V^{i,b}|=1$ and $|V^{i,a}|=|V^{i,a'}|=1$. By Lemma \ref{wa}, $V^{i,l}\sqsubseteq  W^{i,l}$ for $l\in \{a,a'\}$ from where we obtain $|W^{i,l}|\geq 2$ for $l\in \{a,a'\}$ (recall that $W$ contains a twin pair). Since $m(V^{i,b'}_{i^c}\setminus V^{i,b}_{i^c})=6$, by (\ref{mb}), $|W^{i,b'}|=6$. It is easy to check that then there are three words in $W^{i,b'}$ such that no two of them form a twin pair. Thus, by Lemma \ref{po7}$iii$, $|V^{i,b'}|\geq 9$, a contradiction. 

If $k=2$ and the words in $W^{i,b'}$ do not form a twin pair, then, by Lemma \ref{po7}$ii$ and (\ref{mb}), $|V^{i,b'}|\geq 7$ and  $|V^{i,b}|\geq 5$, a contradiction. If $W^{i,b'}$ contains a twin pair, then $|V^{i,b'}|\geq 6$ and  $|V^{i,b}|\geq 4$, a contradiction, because $V^{i,a}\cup V^{i,a'}\neq\emptyset$.

Finally, if $k=1$, then,  by Lemma \ref{po7}$i$ and (\ref{mb}), $|V^{i,b'}|\geq 5$ and  $|V^{i,b}|\geq 4$, and then $V^{i,l}=\emptyset$ for some $l\in \{a,a'\}$ which is not possible.

We showed that $V^{i,l}, W^{i,l}\neq\emptyset$ for every $i\in [d]$ and $l\in S$.

It follows from the above that the distributions $4-8$ are not possible. Indeed, since $K=\breve{u}\cup \breve{w}$ and $u_{1^c}=w_{1^c}$ (because $u,w\in W$ form a twin pair in which $u_1=w'_1$), we have $V^{i,u_i'}=\emptyset$ for $i\in \{2,...,d\}$. 

To show that the distribution 2 is impossible let $U\subset V$ be the code that generates the partition $\ka F$ with the distribution given at the position  2, that is, $U$ consists of all $v\in V$ such that $\breve{v}\cap K\neq\emptyset$ and $|U^{1,a}|=|U^{1,a'}|=2$ and $|U^{1,b}|=|U^{1,b'}|=2$. Consequently, by Lemma \ref{wac}, the codes $U^{1,l}_{1^c}$ and  $U^{1,l'}_{1^c}$ are $u_{1^c}$-equivalent for $l\in \{a,b\}$, and hence they are, up to isomorphism, of the form given in Lemma \ref{cody}$ii$. Thus, for $l\in \{a,a',b,b'\}$ there is $k\in \{0,...,d-3\}$ such that  
\begin{equation}
\label{ppp}
\sum_{v\in U^{1,l}}g(v_{1^c}, u_{1^c})=3\cdot 2^k.
\end{equation}
Since the code $U^{1,a}_{1^c}\cup U^{1,b}_{1^c}$ is a cover of $u_{1^c}$, by (\ref{2d}),  $\sum_{v\in Q}g(v_{1^c}, u_{1^c})=2^{d-1},$ where $Q=U^{1,a}\cup U^{1,b}$.
This, by (\ref{ppp}), is not possible. A contradiction. 

To show that also the rest of the distributions of $\ka F$ is not possible we now prove that $|V^{i,l}|\leq 5$ for every $i\in [d]$ and every $l\in S$. To do this, suppose on the contrary that there are $i\in [d]$ and $l\in S$ such that $|V^{i,l}|\geq 6$. We may assume $l=b$. Since $V^{i,s}\neq\emptyset$ for $s\in \{a,a',b'\}$, we have $m(V^{i,b}_{i^c}\setminus V^{i,b'}_{i^c})\geq 4$. For the same reason we may assume that $|V^{i,a}|=|V^{i,a'}|=1$  or $|V^{i,a}|=1$ and $|V^{i,a'}|=2$, and  then, by Lemma \ref{wa}, $V^{i,s}\sqsubseteq W^{i,s}$ for $s\in \{a,a'\}$. Below we examine these two cases.

Case $|V^{i,a}|=|V^{i,a'}|=1$. As $W$ contains twin pairs, we have $|W^{i,s}|\geq 2$ for $s\in \{a,a'\}$. Since $m(V^{i,b}_{i^c}\setminus V^{i,b'}_{i^c})\geq 4$, it follows that, by (\ref{mb}), $|W^{i,b}|=5$ and $|W^{i,b'}|=1$ or $|W^{i,b}|=4$ and $|W^{i,b'}|=2$.   

Let $|W^{i,b}|=5$ and $|W^{i,b'}|=1$. If there is $v\in V^{i,a}\cup V^{i,a'}$ such that $\breve{v}\cap \bigcup E(W^{i,b}\cup W^{i,b'})\neq\emptyset$, then, by Lemma \ref{po1m}, $W^{i,b'}_{i^c}\sqsubseteq W_{i^c}^{i,b}$, and hence $m(W^{i,b}_{i^c}\setminus W^{i,b'}_{i^c})=4$.
 Since $m(W^{i,b}_{i^c}\setminus W^{i,b'}_{i^c})=m(V^{i,b}_{i^c}\setminus V^{i,b'}_{i^c})$,  $|V^{i,b}|\geq 6$ and $V^{i,l}\neq\emptyset$ for $l\in S$ , we obtain $V_{i^c}^{i,b'}\sqsubseteq V_{i^c}^{i,b}$ and  $|V^{i,b'}|=2$, $|V^{i,b}|=6$. This is not true, by Lemma \ref{po7}$ii$. Thus, $W^{i,b}\sqsubseteq V^{i,b}$ and  $W^{i,b'}\sqsubseteq V^{i,b'}$, which gives, by Lemma \ref{po7}$i$,  $|V^{i,s}|\geq 5$  for $s\in \{b,b'\}$. A contradiction, because $V^{i,l}\neq\emptyset$  for $l\in \{a,a'\}$.

Let  now $|W^{i,b}|=4$ and $|W^{i,b'}|=2$. Since $m(V^{i,b}_{i^c}\setminus V^{i,b'}_{i^c})\geq 4$, by (\ref{mb}), $m(W^{i,b}_{i^c}\setminus W^{i,b'}_{i^c})=4$ and then $m(V^{i,b}_{i^c}\setminus V^{i,b'}_{i^c})=4$. Moreover, $m(W^{i,b'}_{i^c}\setminus W^{i,b}_{i^c})=2$. Consequently, by (\ref{mb}),  $m(V^{i,b}_{i^c}\setminus V^{i,b'}_{i^c})=6$. A contradiction, by (\ref{mb}).

\smallskip
Case $|V^{i,a}|=1$ and $|V^{i,a'}|=2$. In this case we have $|V^{i,b}|=6$ and $|V^{i,b'}|=1$. Moreover,  $|W^{i,a}|\geq 2$ and $|W^{i,a'}|\geq 3$, because, by Lemma \ref{wa}, $V^{i,s}\sqsubseteq W^{i,s}$ for $s\in \{a,a'\}$. Thus, $|W^{i,b}|\leq 4$, because $W^{i,b'}\neq\emptyset$.  Hence, $m(W^{i,b}_{i^c}\setminus W^{i,b'}_{i^c})\leq 4$. This is not possible, because $m(V^{i,b}_{i^c}\setminus V^{i,b'}_{i^c})\geq  5$. 

We proved that $|V^{i,l}|\leq 5$ for every $i\in [d]$ and every $l\in S$. Therefore, the distributions $|\ka F^{1,a}|=m-1$, $|\ka F^{1,a'}|=1$ for $m\in \{7,...,10\}$ are not possible. (Since in the same way we prove that $|V^{i,l}|\leq 5$ for $d=5$, let us notice here that also the distribution $|\ka F^{1,a}|=5$, $|\ka F^{1,a'}|=1$ for $d=5$ is not possible. To explain it, note that, by Lemma \ref{cody}$i$, the code $U=U^{1,a}\cup U^{1,a'}$, where $U^{1,l}=\{v\in V^{1,l}\colon \breve{v}\cap K\neq\emptyset\}$ for $l\in \{a,a'\}$, which generates the set of boxes $\ka F^{1,a}\cup \ka F^{1,a'}$ (and thus, $U^{1,a}_{1^c},U^{1,a'}_{1^c}$ are $u_{1^c}$-equivalent) has to be, up to isomorphism, of the form
\begin{equation}
\label{v5}
U^{1,a}=\{aaaab,aa'a'a'b,abaa'b,aa'bab,aaa'bb\},\;\; U^{1,a'}=\{a'bbbb\},
\end{equation}
where we assumed that $u_{1^c}=bbbb$. Then, $|U^{5,b}|=6$, a contradiction.) 

\smallskip
Skipping the last letter $b$ in (\ref{v5}) we obtain the code $U=U^{1,a}\cup U^{1,a'}$ which generates the partition $\ka F=\ka F^{1,a}\cup \ka F^{1,a}$ with $|\ka F^{1,a}|=5$, $|\ka F^{1,a'}|=1$ for $d=4$. Computing all twin pair free codes  $V$ such that $U\subset V$, $|V|=10$ and $|V^{i,l}|\leq 5$ for $i\in [4]$ and $l\in S$ we always obtain   $V^{i,l}=\emptyset$ for some $i\in [4]$ and some $l\in S$. Thus, $\ka F$ cannot have the distribution 1 for $m=6$.

If $\ka F$ has the distribution number 3, then there is $U\subset V$ such that $|U^{1,b}|=2,|U^{1,b'}|=3$ and $U$ generates $\ka F^{1,b}\cup \ka F^{1,b'}$ (this case is discussed for $d=4,5$). Since $U^{1,b}_{1^c}$ and $U^{1,b}_{1^c}$ are $u_{1^c}$-equivalent, by Lemma \ref{cody}$iii$,  

\vspace{-2mm}
\begin{equation}
\label{v231}
U^{1,b}=\{babbs,\; ba'a'a's\},\;\; U^{1,b'}=\{b'ba'a's,\; b'abas,\; b'aaa's\},
\end{equation}

\vspace{-4mm}
or

\vspace{-4mm}
\begin{equation}
\label{v232}
U^{1,b}=\{baabs,\; bba'a's\},\;\; U^{1,b'}=\{b'aba's,\; b'aaas,\; b'a'a'a's\},
\end{equation}
where $s\in \{a,a',b\}$ (it is assumed that $u_{1^c}=b...b\in S^{d-1}, d=4,5$). Now we compute all twin pair free codes $Q$ such that $\ka F=\{K\cap \breve{v}\colon v\in Q\}$, that is, $U\subset Q$, $|Q|=9$, $|Q^{1,a}|=|Q^{1,a'}|=2$, and moreover the codes $Q^{1,a}_{1^c}\cup Q^{1,b}_{1^c}$ and $Q^{1,a'}_{1^c}\cup Q^{1,b'}_{1^c}$ are  covers of the word $u_{1^c}=bbbb$ if $d=5$ and $u_{1^c}=bbb$ if $d=4$ (for $d=4$ the codes $U$ are obtained from (\ref{v231}) and (\ref{v232}) by skipping the letter $s$). The computations show that there is, up to isomorphism, one such code (we give the form for $d=5$; skipping the last letter in $Q$ we get the form for $d=4$):
$$
Q=\{aa'aa'b, aa'bab, a'a'abb, a'a'a'ab, babbb, ba'a'a'b, b'ba'a'b, b'abab, b'aaa'b\},
$$
but  for every $v\in S^d$ if  $V=Q\cup \{v\}$ is a twin pair free code, then $V^{i,l}=\emptyset$ for some $i\in [d]$ and $l\in S$. A contradiction.
\hfill{$\square$}

\smallskip
Recall that a pair of dichotomous words $u,v$ is an $i$-siblings if $u,v$ are not a twin pair but $u_{i^c},v_{i^c}$ form a twin pair. 
\begin{wn}
\label{sy1}
Let $V,W\subset S^d$, where $S=\{a,a',b,b'\}$ and $d=5,6$, be disjoint equivalent polybox codes without twin pairs, and let $|V|\in \{13,...,16\}$. For every $i\in [d]$ and every $l,s\in S$, $l\not\in \{s,s'\}$, if  $4\leq |V^{i,l}\cup V^{i,s}|\leq 10$, then the code  $V^{i,l}\cup V^{i,s}$ contains an $i$-siblings.

\end{wn}
\proof Suppose on the contrary that there are no $i$-siblings in $V^{i,l}\cup V^{i,s}$ for some $i\in [d],l,s\in S$, $l\not\in \{s,s'\}$. By Lemma \ref{r24} and \ref{rr2}, the sets $V^{i,l},V^{i,s}$ are non-empty. By Lemma \ref{rig}, the code $V_{i^c}^{i,l}\cup V_{i^c}^{i,s}$ is rigid, and thus, by Lemma \ref{cy}, $V_{i^c}^{i,l}\cup V_{i^c}^{i,s}=W_{i^c}^{i,l}\cup W_{i^c}^{i,s}$. Since $V\cap W=\emptyset$, we have $V_{i^c}^{i,l}=W_{i^c}^{i,s}$ and $V_{i^c}^{i,s}=W_{i^c}^{i,l}$. Then, by ({\bf Co}),  $W_{i^c}^{i,s}\sqsubseteq V_{i^c}^{i,l'}$ and $W_{i^c}^{i,l}\sqsubseteq V_{i^c}^{i,s'}$. Since $|V^{i,l}\cup V^{i,s}|\geq 4$, by Lemma \ref{po7}, $|V^{i,l'}\cup V^{i,s'}|\geq 14$, a contradiction. 
\hfill{$\square$}

\begin{wn}
\label{le12}
Let $V,W\subset S^d$, where $S=\{a,a',b,b'\}$ and $d=5,6$, be disjoint equivalent polybox codes without twin pairs, $|V|\in \{13,...,16\}$, and let $V$ be not flat.
If $d=6$, then $|\{v\in V:\breve{v}\cap \breve{w}\neq \emptyset\}|\leq 10$ for every $w\in W$, and if $d=5$, then $|\{v\in V:\breve{v}\cap \breve{w}\neq \emptyset\}|\leq 11$ for every $w\in W$.
\end{wn}
\proof Let $d=6$ and $V_w=\{v\in V:\breve{v}\cap \breve{w}\neq \emptyset\}$.  
Suppose on the contrary that $|V_w|\geq 11$ for some $w\in W$. We assume that $|V|=16$ and $|V_w|=11$. Without loss of generality we may take $w=b...b$. Then every word from the set $V\setminus V_w=\{v^1,v^2,v^3,v^4,v^5\}$ contains the letter $b'$ at some position $i\in [d]$ and there is no word with the letter $b'$ in the set $V_w$.

Now we consider the graph of siblings $G=(V,\ka E)$.
In what follows we shall use the notation used in the proof of Lemma 22 of \cite{Kis22}: An edge between vetrices $v^i$ and $v^j$, $i,j\in [5]$, is called {\it internal}, while an edge joining $v\in V\setminus V_w$ and $u\in V_w$ is called {\it external}. An edge between $v,u\in V$ is of the {\it type $b'$ with the colour $i$}, if $v_i=b'$ and $u_i\in \{a,a'\}$ or $u_i=b'$ and $v_i\in \{a,a'\}$ for some $i\in [d]$. Such an edge is denoted by $(v,u)_{b'}$.   


Note that each edge of the type $b'$ has to be incident to a vertex from the set $V\setminus V_w$, as there is no word with the letter $b'$ in the set $V_w$. It is easy to check that if $v\in V\setminus V_w$ has more than two letters $b'$, then there is no external edge which is incident to $v$; if $v$ has exactly two letters $b'$, then there is at most one external edge incident to $v$. Finally, if $v\in V\setminus V_w$ has one letter $b'$, then all external edges of the type $b'$ which are incident to $v$ are of the same colour. 
 
Recall that, by Lemma \ref{r24} and \ref{rr2}, $V^{i,l}\neq\emptyset$ for every $i\in [d]$ and $l\in S$.

We show that for every $i\in [d]$ there is an edge of the type $b'$ in $G$ with the color $i$. To do this, we consider two cases: $|V^{i,l}\cup V^{i,b'}|\leq 3$ for some $l\in \{a,a'\}$ and $|V^{i,l}\cup V^{i,b'}|\geq 4$ for $l\in \{a,a'\}$.

Let $|V^{i,l}|=1, |V^{i,b'}|=2$. If there are no $i$-siblings in the set $V^{i,l}\cup V^{i,b'}$, then we show, in the same way as in the proof of Lemma \ref{sy1}, that $|V^{i,l'}|\geq 5$ and $|V^{i,b}|\geq 7$. The code $V$ has sixteen words, and therefore $|V^{i,l'}|\leq 6$. Then, by Lemma \ref{sy1}, the code $V^{i,l'}\cup V^{i,b'}$ contains an $i$-siblings.
Similarly we show that  the set $V^{i,l}\cup V^{i,b'}$ or $V^{i,l'}\cup V^{i,b'}$ has to contain an $i$-siblings if $|V^{i,l}|=2$ and $|V^{i,b'}|=1$ or $|V^{i,l}|=|V^{i,b'}|=1$ for some $l\in \{a,a'\}$.

Let now $|V^{i,l}\cup V^{i,b'}|\geq 4$ for $l\in \{a,a'\}$. Since $1\leq |V^{i,b'}|\leq 5$ and $|V|=16$, we have  $4\leq |V^{i,a}\cup V^{i,b'}|\leq 10$ or $4\leq |V^{i,a'}\cup V^{i,b'}|\leq 10$. Consequently, by Lemma \ref{sy1}, there is an $i$-siblings in $V^{i,a}\cup V^{i,b'}$ or in $V^{i,a'}\cup V^{i,b'}$.
 
Thus, for every $i\in [d]$ there is an edge of the type $b'$ in $G$ with the colour $i$.  

It is easy to check that for every $i,j\in [5],i\neq j$, and every $v,u\in V_w, v\neq u$, the three edges $(v,v^i)_{b'},(v^i,v^j)_{b'}$ and $(v^j,u)_{b'}$ cannot have three different colours. Hence, if every two internal edges are not incident, then there is $i\in [d]$ such that there is no edge of the type $b'$ in $G$ with the colour $i$ (recall that $V$ is a twin pair free code). This, as we showed,  cannot happen and therefore we may assume that there are at least two incident internal edges of the type $b'$. Let $(v^1,v^2)_{b'}$ and $(v^2,v^3)_{b'}$ be these edges.   Note that we may assume that $v^1\in \{b'aaaaa,...,b'b'b'b'b'b'\}$, and then we can enumerate all, up to isomorphism (of polybox codes), edges $(v^1,v^2)_{b'}$. Having this, for every  $(v^1,v^2)_{b'}$ we can compute all edges  $(v^2,v^3)_{b'}$.  (For every $v^1\in  \{b'aaaaa,...,b'b'b'b'b'b'\}$ the number of twin pair free codes $\{v^1,v^2,v^3\}$ with the above properties ranges between 36 and 84.)

Now we consider two cases:

1. Let $A=\{v^1,v^2,v^3\}$ and $B=\{v^4,v^5\}$ and suppose that every two vertices $u\in A$ and $v\in B$ are not joined by an edge of the type $b'$. Since the total number of edges of the type $b'$ with different colours which are incident to a vertex from $B$ is less than three, to obtain six edges of the type $b'$ in all six colours, we need at least four edges of the type $b'$ with four different colours which are incident to vertices from $A$. The computations show that this is impossible.

2. Now we assume that  there is no edge of the type $b'$ between every two vertices $u,v$ from the sets $A=\{v^1,v^2,v^3,v^4\}$ and $B=\{v^5\}$, respectively, and moreover every two vertices in $A$ are connected by a path consisting of edges of the type $b'$. Now we need at least five edges of the type $b'$ with five different colours which are incident to vertices from $A$. Also in this case the computations show that it is not possible.

It follows from the above that every two vertices in the set $\{v^1,...,v^5\}$ are connected by a path consisting of edges of the type $b'$. 
The number of the codes $\{v^1,v^2,v^3\}$ is low, and thus it is easy to compute all such sets $\{v^1,...,v^5\}$, which are, recall that, twin pair free codes. Since $V^{i,l}\neq\emptyset$ for $i\in [d]$ and $l\in S$, for every code $\{v^1,...,v^5\}$ and every $i\in [d]$, there is $j\in [5]$ such that $v^j_i=b'$.
(Note that we can eliminate from the further computations all codes  $\{v^1,...,v^5\}$ for which the sum of the number of internal edges of the type $b'$ with different colours and the number of words in $\{v^1,...,v^5\}$ which have less then three letters $b'$ each, is at most five.)
Let $e(v^i)\subset V_w$ consist of all $v\in V_w$ that form with $v^i$ an edge of the type $b'$ and $\{v^1,...,v^5,v\}$ is a twin pair free code. Every set  $e(v^i)$, $i\in [5]$, has a few words, and therefore for every $\{v^1,...,v^5\}$ we can compute twin pair free codes  $\{v^1,...,v^5\}\cup P$, where $P\subset \bigcup_{i=1}^5e(v^i)$. In each computed  code $\{v^1,...,v^5\}\cup P$ lacks an edge of the type $b'$ with some colour $i\in [6]$. A contradiction. 

Clearly, we also obtain a contradiction if $|V|=16$ and  $|V_w|\geq 12$  or $|V|\leq 15$ and  $|V_w|\geq 11$

Along the same lines we prove the lemma for $d=5$.
\hfill{$\square$}

\section{The number of letters}

In this section we show that the codes $V$ and $W$ under consideration can be written down in the alphabet $S=\{a,a',b,b'\}$.  We begin with the following lemma:

\begin{lemat}
\label{ogr}
Let $V,W\subset S^d$, where $S=\{a,a',b,b',c,c'\}$ and $d=5,6$, be disjoint equivalent polybox codes without twin pairs, $|V|\in \{13,...,16\}$, and let $V$ be  not flat. Assume that there is $i \in [d]$ such that $V^{i,l}\cup V^{i,l'}\neq\emptyset$ for every $l\in \{a,b,c\}$. Then there is $l\in S$ such that $|V^{i,l}|=|V^{i,l'}|=2$ and $|W^{i,l}|,|W^{i,l'}|\leq 3$ (or conversely with respect to $V$ and $W$) or $|V^{i,l}|=2$, $|V^{i,l'}|=3$ and $|W^{i,l}|=2, |W^{i,l'}|=3$. Moreover, every word in $V^{i,l}\cup V^{i,l'}$ (or in $W^{i,l}\cup W^{i,l'}$) is covered by at most eleven words from $W$ (or from $V$). 
\end{lemat} 
\proof By Lemma \ref{rr2}, ({\bf C}) and Theorem \ref{t12}, the sets $V^{i,l}$ and $W^{i,l}$ are non-empty for every $l\in S$.

Suppose that for some $l\in S$, say $l=a$, we have $|V^{i,a}|=|V^{i,a'}|=1$. Then, by Lemma \ref{wa} and \ref{po7}$i$, $|W^{i,a}|\geq 5$ and $|W^{i,a'}|\geq 5$. Since $W^{i,l}\neq\emptyset$ for $l\in \{b,b',c,c'\}$, we may assume, by Lemma \ref{rr3}$i$, that $|W^{i,a}|=|W^{i,a'}|=5$, $|W^{i,b}|=|W^{i,b'}|=2$ and $|W^{i,c}|=|W^{i,c'}|=1$. Then  $|V^{i,c}|=|V^{i,c'}|=5$ and $|V^{i,b}|=|V^{i,b'}|=2$.

Since the distribution $|V^{i,l}|=1$ and $2\leq |V^{i,l'}|\leq 4$, by Lemma \ref{rr3}$i$, is not possible, if $4\leq |V^{i,l}\cup V^{i,l'}|\leq 5$, then  $|V^{i,l}|=|V^{i,l'}|=2$ or $|V^{i,l}|=2$, $|V^{i,l'}|=3$.  Clearly,  $|V^{i,l}\cup V^{i,l'}|\leq 5$ for at least one $l\in \{a,b,c\}$, and thus $|V^{i,l}|=|V^{i,l'}|=2$ or $|V^{i,l}|=2$, $|V^{i,l'}|=3$. We now consider these two cases assuming that $l=a$ and $i=1$. 

Let $d=6$.

If $|V^{1,a}|=|V^{1,a'}|=2$, then, by Lemma \ref{rr3}$ii$, there is $w\in W\setminus (W^{1,a}\cup W^{1,a'})$ such that $\breve{w}\cap \bigcup E(V^{1,a}\cup V^{1,a'})\neq\emptyset$. Thus, by Lemma \ref{wa1} and  \ref{wac}, the codes $V_{1^c}^{1,a}, V_{1^c}^{1,a'}$ are $w_{1^c}$-equivalent, and consequently, by Lemma \ref{cody}$ii$, up to isomorphism,
\begin{equation}
\label{v2}
V^{1,a}\cup V^{1,a'}=\{aabbbb,\; aa'abbb,\;a'babbb,\;a'aa'bbb\}.
\end{equation}
Observe that, $m(V_{1^c}^{1,l}\setminus V_{1^c}^{1,l'})=\frac{3}{4}$ for $l\in \{a,a'\}$ and  therefore, by (\ref{mb}), $|W^{1,a}|=|W^{1,a'}|=m$ for some $m\geq 1$. By the previous case (that is, $|V^{1,a}|=|V^{1,a'}|=1$) we assume $m\leq 4$.  If $m=4$, then $|W^{1,l}|=|W^{1,l'}|=2$ for $l\in \{b,c\}$, by Lemma   \ref{rr3}$i$ (compare the case $|V^{1,a}|=|V^{1,a'}|=1$ to see that it cannot be $|W^{1,l}|=|W^{1,l'}|=1$ for some $l\in \{b,c\}$). Thus, again  by (\ref{mb}), $|V^{1,l}|=|V^{1,l'}|$ for $l\in \{b,c\}$. Since $|V|\leq 16$, we have $|V^{1,l}|,|V^{1,l'}|\leq 3$ for $l=b$ or $l=c$. Thus, if $m=4$, then $|W^{i,l}|=|W^{i,l'}|=2$ and $|V^{i,l}|,|V^{i,l'}|\leq 3$ for $l=b$ or $l=c$, and if $m\leq 3$, then $|V^{i,a}|=|V^{i,a'}|=2$ and $|W^{i,a}|,|W^{i,a'}|\leq 3$.


\smallskip
Suppose now that $|V^{1,a}|=2$ and $|V^{1,a'}|=3$. It follows from the previous case that we may assume that $|V^{1,l}\cup V^{1,l'}|\geq 5$ and $|W^{1,l}\cup W^{1,l'}|\geq 5$ for every $l\in S$. 
Since, by Lemma \ref{rr3}$i$,   $|V^{1,l}|=2$, $|V^{1,l'}|=3$ for at least one $l\in \{b,c\}$ and $|W^{1,l}|=2$, $|W^{1,l'}|=3$ for at least two letters $l\in \{a,b,c\}$, the proof of the first part of the lemma is completed.

\smallskip
To prove the second part of the lemma assume first that $|V^{1,a}|=|V^{1,a'}|=2$ and $|W^{1,l}|\leq 3$ for $l\in \{a,a'\}$. Clearly, the code $V^{1,a}\cup V^{1,a'}$ is of the form (\ref{v2}). Let $v=aa'abbb$ and $u=a'babbb$. We shall consider covers $C_u\subset W$ and $C_v\subset W$ of $u$ and $v$, respectively.

{\center
\includegraphics[width=5cm]{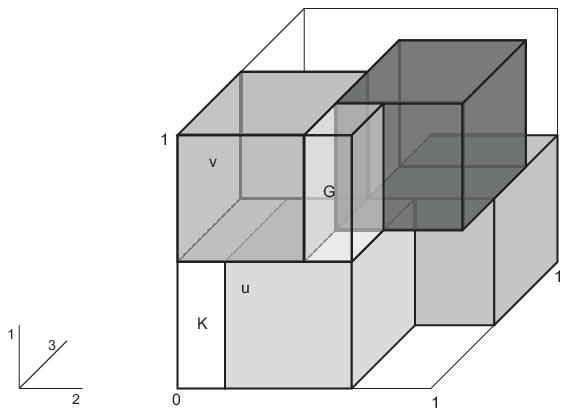}\\
}

\medskip
\noindent{\footnotesize 
Figure 5: A scheme of the realization $E(V^{1,a}_H\cup V^{2,a'}_H)$, where $H=\{1,2,3\}$. In the picture only the parts $K_H$ and $G_H$ of the boxes $K$ and $G$ are visible. In this scheme we identify $Ea$ with the interval $[\frac{1}{2},1]$ on the first and second axis and with $[0,\frac{1}{2})$ on the third axis; 
$Eb$ is identified with the interval $[\frac{1}{4},\frac{3}{4}]$ on all axes.
  }

\smallskip
Assume on the contrary that the word $u$ is covered by at least twelve words from $W$. Note that, by ({\bf P}), the box $K=Ea'\times (Ea'\cap Eb')\times Ea\times (Eb)^3$ (Figure 5) cannot be intersected by a box $\breve{w}$ for  $w\in V\cup W$. Therefore, $w_2=b$ for every word $w\in C_u\subseteq W$. Since $|C_u|\geq 12$, we obtain a contradiction, by Lemma \ref{rr3}$i$.

Now we assume that $v$ is covered by a code $C_v\subset W$ having at least twelve words. By Lemma \ref{rr3}$ii$, the set $P=\{w\in  W\setminus (W^{1,a}\cup W^{1,a'})\colon \breve{w}\cap \bigcup E(V^{1,a}\cup V^{1,a'})\neq\emptyset\}$ is non-empty. 
Moreover, by Lemma \ref{cody}$iv$, $w_2=w_3=b$ for every $w\in P$. 
Since, $|W^{1,l}|\leq 3$ for $l\in \{a,a'\}$ and $P\subseteq C_v$, we have $|P|\geq 9$.  By Lemma \ref{ci} (in which we take $A=\{2,3\}, B=\{4,5,6\}, p=bb, q=bbb$ and, by Lemma \ref{wa1}, $U^{1,l}=V^{1,l}$ for $l \in \{a,a'\}$), the codes $P_{B}^{1,b}\cup P_{B}^{1,c}$ and $P_{B}^{1,b'}\cup P_{B}^{1,c'}$ are $q$-equivalent. 
Let $K,M$ be such as in the first row of Table 2, $K^1=\{a'bb,aa'a'\}$, $M^1=\{ba'a',a'ab,a'a'a\}$, $K^2=\{ba'a',a'ba\}$, $M^2=\{a'a'b,a'aa,aa'a'\}$ and $K^3, M^3$ be such as in the last row of Table 2 (the codes $K^1,M^1$ and $K^2,M^2$ are isomorphic to $K,M$ in the third and the fourth row of Table 2, respectively).

Since $|P|\geq 9$, an inspection of the codes in Table 2 show that, up to isomorphism, there are three possibilities:

1. $P_{B}^{1,b}\cup P_{B}^{1,b'}=K\cup M$ and $P_{B}^{1,c}\cup P_{B}^{1,c'}=K^1\cup M^1$.

2. $P_{B}^{1,b}\cup P_{B}^{1,b'}=K\cup M$ and $P_{B}^{1,c}\cup P_{B}^{1,c'}=K^2\cup M^2$.

3. $P_{B}^{1,b}\cup P_{B}^{1,b'}=K^3\cup M^3$ and $P_{B}^{1,c}\cup P_{B}^{1,c'}=\emptyset$. 

\smallskip
Now we compute all twin pair free covers $C_v$ of $v$ containing $P$. Observe that, by the condition $|W^{1,l}|\leq 3$ for $l\in \{a,a'\}$ we have $|C_v|=12$ in the cases 1-2; if $P$ is such as in the third case, we have $|C_v|=12$ or $|C_v|=13$. By (${\bf P}$), every such $C_v$ cannot contain a word $w$ with $\breve{w}\cap K\neq\emptyset$ or $\breve{w}\cap G\neq\emptyset$, where  $G=Ea\times (Ea\cap Eb)\times (Ea\cap Eb')\times (Eb)^3$ (Figure 5). Easy computations (which can be made even by hand) show that every computed cover $C_v$ with these properties is of the form $C=C^{3,b}\cup C^{3,b'}$. Since $C_v\subset  W^{3,b}\cup W^{3,b'}$, we have $|W^{3,b}\cup W^{3,b'}|\geq 12$. But note that for every $w\in P$ the box $\breve{w}$ intersects $\breve{u}$, that is, $P\subset C_u$ for every cover $C_u\subset W$ of $u$. For every $P$ mentioned in 1-3 we compute twin pair free cover $C_u$ such that $|C_u|\leq 11$. In every such $C_u$ there is a word $w$ with $w\in W^{3,b'}$ and $w_1=a'$.  Thus, $w\not\in C_v$ for every above computed cover $C_v$ of $v$. Therefore, $|W^{3,b}\cup W^{3,b'}|\geq 13$. Then, by Lemma \ref{rr3}$i$, $|W^{3,a}|=|W^{3,a'}|=1$. This, by Lemma \ref{rr5}, is not possible, because $P\subseteq W^{3,b}$ and $|P|\geq 9$.



\medskip
Let now $|V^{1,a}|=2$ and $|V^{1,a'}|=3$ and $|W^{1,l}|\leq 3$ for $l\in \{a,a'\}$. 
We shall show, by Lemma \ref{cody}$iii$, that, up to isomorphism, 

\vspace{-2mm}
\begin{equation}
\label{va}
V^{1,a}=\{aabbbb,\; aa'a'a'bb\},\;\; V^{1,a'}=\{a'ba'a'bb,\; a'ababb,\; a'aaa'bb\},
\end{equation}

\vspace{-4mm}
or

\vspace{-4mm}
\begin{equation}
\label{vb}
V^{1,a}=\{aaabbb,\; aba'a'bb\},\;\; V^{1,a'}=\{a'aba'bb,\; a'aaabb,\; a'a'a'a'bb\},
\end{equation}
It follows from Lemma \ref{rr3}$ii$ that there is a word $w\in   W^{1,l}\cup W^{1,l'}$, where $l\in \{b,c\}$, such that the sets  $U_w^{1,a}=\{v\in  V^{1,a}\colon \breve{w}\cap \breve{v}\neq\emptyset\}$ and $U_w^{1,a'}=\{v\in  V^{1,a'}\colon \breve{w}\cap \breve{v}\neq\emptyset\}$ are non-empty, and by Lemma \ref{wac} and \ref{wa1}, $|U_w^{1,a}|\geq 2, |U_w^{1,a'}|\geq 2$. To prove that the codes $V^{1,a}, V^{1,a'}$ have the form (\ref{va}) or (\ref{vb}), it is enough to exclude the case  $|U_w^{1,a}|=|U_w^{1,a'}|=2$. To do this, we assume on the contrary that $|U_w^{1,a}|=|U_w^{1,a'}|=2$ for every $w\in W$ for which the sets $U_w^{1,a},U_w^{1,a'}$ are non-empty. Clearly, by Lemma \ref{wac} and \ref{cody}$ii$, up to isomorphism,
$$
U_w^{1,a}=\{aabbbb,aa'abbb\},\;\; U_w^{1,a'}=\{a'babbb,a'aa'bbb\}.
$$
Let $\{v\}=V^{1,a'}\setminus U_w^{1,a'}$. If $v\in U^{1,a'}_{\bar{w}}$ for some $\bar{w}\in W^{1,l}\cup W^{1,l'}$, $\bar{w}\neq w$, $l\in \{b,c\}$, then $v=a'l_1l_2bbb$ (because $U^{1,a}_{\bar{w}}=U^{1,a}_{w}=V^{1,a}$). It is easy to see that there are no letters $l_1,l_2\in S$ such that $U^{1,a'}_w\cup \{v\}$ is a twin pair free code. Therefore, $v\sqsubseteq W^{1,a'}$, and consequently, $|W^{1,a'}|\geq 5$, which contradicts the assumption $|W^{1,a'}|\leq 3$. Thus, $V^{1,a},V^{1,a'}$ are of the form (\ref{va}) or (\ref{vb}).

\smallskip
Let  $P=\{w\in  W\setminus (W^{1,a}\cup W^{1,a'})\colon \breve{w}\cap \bigcup E(V^{1,a}\cup V^{1,a'})\neq\emptyset\}$. Similarly like above, $P\neq \emptyset$ and, by Lemma \ref{cody}$iv$, $w_2=w_3=w_4=b$ for every $w\in P$. Since now, by Lemma \ref{ci}, the codes $P_{B}^{1,b}\cup P_{B}^{1,c}$ and $P_{B}^{1,b'}\cup P_{B}^{1,c'}$ are $q$-equivalent, where $B=\{5,6\}$ and $q=bb$, we have  $|P|\leq 4$ (this means, that only one of the sets $P_{B}^{1,b}\cup P_{B}^{1,b'}$ and $P_{B}^{1,c}\cup P_{B}^{1,c'}$ is non-empty). Therefore, in every cover $C\subseteq W$ of a word  $v\in V^{1,a}\cup V^{1,a'}$ there are at most four words with the first letter different form $a$ or $a'$. Consequently, if $|C|\geq 12$, then $|C^{1,a}|\geq 4$ or $|C^{1,a'}|\geq 4$. This contradicts the assumption  $|W^{1,a}|,|W^{1,a'}|\leq 3$. 

If $d=5$ and $|V^{1,a}|=|V^{1,a'}|=2$, where $|W^{1,a}|,|W^{1,a'}|\leq 3$, then the codes $V^{1,a},V^{1,a'}$ are obtained from (\ref{v2}) by skipping the last letter $b$. Since now $w_2=w_3=b$ for every $w\in P$, where $P=\{w\in  W\setminus (W^{1,a}\cup W^{1,a'})\colon \breve{w}\cap \bigcup E(V^{1,a}\cup V^{1,a'})\neq\emptyset\}$, similarly like above we obtain $|C^{1,a}|\geq 4$ or $|C^{1,a'}|\geq 4$ in every cover $C\subseteq W$ of a word  $v\in V^{1,a}\cup V^{1,a'}$ with $|C|\geq 12$. A contradiction.

Skipping the last letter $b$ in every word in (\ref{va}) and (\ref{vb}) we obtain the codes $V^{1,a},V^{1,a'}$ in the case  $|V^{1,a}|=2$, $|V^{1,a'}|=3$ for $d=5$. Since now $w_2=w_3=w_4=b$ for every $w\in P$, where $P=\{w\in  W\setminus (W^{1,a}\cup W^{1,a'})\colon \breve{w}\cap \bigcup E(V^{1,a}\cup V^{1,a'})\neq\emptyset\}$, by Lemma \ref{ci} (in which we take $A=\{2,3,4\}$, $B=\{5\}$, $p=bbb$ and $q=b$), the code $V$ or $W$ contains a twin pair, which is not true. The proof of the lemma is completed.
\hfill{$\square$}

\medskip
We are ready to prove the announced restriction on $S$: 
\begin{lemat}
\label{3li}
Let $V,W\subset S^d$, where $d=5,6$, be disjoint equivalent polybox codes without twin pairs, $|V|\in \{13,...,16\}$, and let $V$ be not flat. Then the codes $V,W$ can be written down in the alphabet $S=\{a,a',b,b'\}$.
\end{lemat}

\proof 
Suppose that the lemma is not true. We shall consider two cases: $|S|\geq 8$ and $|S|=6$.

\smallskip
{\bf Case $|S|\geq 8$}. Assume that  there is $i\in [d]$ such that $V^{i,l}\cup V^{i,l'}\neq\emptyset$ for every $l\in S_1=\{a,a',b,b',c,c',d,d'\}$, where $S_1\subseteq S$. Clearly, as we have seen before,  $W^{i,l}\cup W^{i,l'}\neq\emptyset$ for every $l\in S_1$.
We shall show that then  $V^{i,l}\neq\emptyset$ and $V^{i,l'}\neq\emptyset$ for $l\in S_1$. 
To do this, we may assume on the contrary that $V^{i,a'}=\emptyset$. Then, $V^{i,a}\sqsubseteq W^{i,a}$. Note that if $|V^{i,a}|=1$, then by Lemma \ref{po7}$i$ and (\ref{mb}), $|W^{i,a}|\geq 5$ and $|W^{i,a'}|\geq 4$ and if $|V^{i,a}|\geq 2$, then,  by Lemma \ref{po7}$ii$ and (\ref{mb}), $|W^{i,a}\cup W^{i,a'}|\geq 8$. Therefore, we may assume that  $|W^{i,b}\cup W^{i,b'}|\leq 3$ and $|W^{i,c}\cup W^{i,c'}|\leq 3$. Now it is easy to show that  $|V^{i,b}\cup V^{i,b'}|\geq 8$ and $|V^{i,c}\cup V^{i,c'}|\geq 8$. A contradiction. Thus,  $V^{i,l}\neq\emptyset$ and $V^{i,l'}\neq\emptyset$ for $l\in S_1$. Clearly, we have also $W^{i,l}\neq\emptyset$ and $W^{i,l'}\neq\emptyset$ for $l\in S_1$.


Note that $|V^{i,l'}\cup V^{i,l}|=4$ for $l \in S_1$. Indeed, if $|V^{i,l'}\cup V^{i,l}|<4$ for some $l\in S_1$, then, by Lemma  \ref{wa} and \ref{po7}$i$,  $|W^{i,l}|, |W^{i,l'}|\geq 5$ and then $W^{i,s}=\emptyset$ for some $s\in S_1$, which cannot occur. 

 
 Consequently,  $|V^{i,l'}\cup V^{i,l}|=4$ for every $l\in S_1$, and then, again by Lemma \ref{wa} and \ref{po7}$i$,  $|V^{i,l'}|=|V^{i,l}|=2$ for $l\in S_1$. Similarly, $|W^{i,l'}|=|W^{i,l}|=2$ for $l\in S_1$. 
Thus, for every $l\in S_1$ there is $w\in W$ with $w_i\not\in \{l,l'\}$  such that $\breve{w}\cap \bigcup E(V^{i,l}\cup V^{i,l'})\neq\emptyset$ (as if it is not so for some $l\in S_1$, then $V^{i,l}\sqsubseteq W^{i,l}$ and  consequently, by Theorem \ref{t12}, $|W^{i,l}|\geq 3$. A contradiction.)
This means, by Lemma \ref{wac}, that the codes $V^{i,l}_{i^c}$, $V^{i,l'}_{i^c}$ are $w_{i^c}$-equivalent. In what follows we may assume that $i=1$, $l=a$ and let $d=6$ (along the same lines we prove the lemma for $d=5$). 
By Lemma \ref{cody}$ii$, we may assume that the code $V^{1,a}\cup V^{1,a'}$ is of the form (\ref{v2}). Then, by Lemma \ref{cody}$iv$,  $w_2=w_3=b$ for every $w\in P$, where  $P=\{w\in  W\setminus (W^{1,a}\cup W^{1,a'})\colon \breve{w}\cap \bigcup E(V^{1,a}\cup V^{1,a'})\neq\emptyset\}$. For every $l\in S\setminus \{a,a'\}$ the codes $P^{1,l}_B, P^{1,l'}_B$, by Lemma \ref{ci}, are $q$-equivalent for $q=bbb$, and then $|P^{1,l}_B|=|P^{1,l'}_B|=2$, where $B=\{4,5,6\}$, because $|W^{i,l}|=2$ for $l\in S_1$. In the similar way as in the middle part of the proof of Lemma \ref{rig} (that part where (\ref{ppp}) is considered) we show that the code $P^{1,b}_B\cup P^{1,c}_B\cup P^{1,d}_B$ cannot cover the word $q$.
Let $(t_4,t_5,t_6)\in \breve{q}\setminus \bigcup \{\breve{w}_B\colon w\in P\}$. Then, the boxes $\breve{v}^1_H\times \{t_4\}\times \{t_5\}\times \{t_6\}$  and $\breve{v}^2_H\times \{t_4\}\times \{t_5\}\times \{t_6\}$, where $\{v^1,v^2\}=V^{1,a}$ and $H=\{1,2,3\}$, must be covered by the boxes $\breve{w}^1,\breve{w}^2$, where $\{w^1,w^2\}=W^{1,a}$. Clearly, it can be done only if $v^i_H=w^i_H$ for $i\in \{1,2\}$, because $|W^{i,a}|=2$. Since the structure of the code $W^{1,a}\cup W^{1,a'}$ is, up to isomorphism, the same as $V^{1,a}\cup V^{1,a'}$, we have $w^i_B=r$ for $i\in \{1,2\}$, where $r\in S^3$. (Obviously, $q\neq r$, because $V$ and $W$ are disjoint.) Note that $\breve{r}\cap \breve{w}_B=\emptyset$ for every $w\in P$, otherwise $\breve{w}^i\cap \breve{w}\neq\emptyset$ for some $w^i\in W^{1,a}$ and some  $w\in P$, which is not possible. Take $(x_4,x_5,x_6)\in \breve{w}_B\cap \breve{q}$ for any fixed $w\in P$ and, by ({\bf P}), $(x_1,x_2,x_3)\in \breve{v}^1_H\setminus \bigcup \{\breve{w}_H\colon w\in  W\setminus (W^{1,a}\cup W^{1,a'})\}$. The point $x=(x_1,...,x_6)$ belongs to $\breve{v}^1$, because $v^1_B=q$. On the other hand $x\not \in \breve{w}^1\cup \breve{w}^2$, because $(x_4,x_5,x_6)\in \breve{w}_B$, where $w\in P$, and $\breve{r}\cap \breve{w}_B=\emptyset$ for every $w\in P$. Clearly, $x\not \in \bigcup E(W^{1,a'})$, because $x_1\in Ea$. Moreover, by the manner of choosing of $(x_1,x_2,x_3)$,  we have $x \not\in \bigcup E( W\setminus (W^{1,a}\cup W^{1,a'}))$. Thus, $x\not \in \bigcup E(W)$, a contradiction. 
This proves that $|S|<8$, that is, $|S|\leq 6$. 

\smallskip
{\bf Case $|S|=6$}. Let $S=\{a,a',b,b',c,c'\}$, and let $V^{1,l}\cup V^{1,l'}\neq\emptyset$ for $l\in S$. Clearly, by ({\bf C}) and Theorem \ref{t12},  $W^{1,l}\cup W^{1,l'}\neq\emptyset$ for $l\in S$, and consequently, by Lemma \ref{rr2}, $V^{1,l},W^{1,l}\neq\emptyset$ for $l\in S$.  By Lemma \ref{ogr} we may assume that $|V^{1,a}|=|V^{1,a'}|=2$ or $|V^{1,a}|=2, |V^{1,a'}|=3$ and, in both cases,  $|W^{1,a}|,|W^{1,a'}|\leq 3$. In the proof of that lemma we showed that, up to isomorphism, the code $V^{1,a}\cup V^{1,a'}$ is of the form (\ref{v2}), (\ref{va}) or (\ref{vb}). In what follows we compute covers of $V^{1,a}\cup V^{1,a'}$ by words from $W$. Clearly, we may assume that $V^{1,a}\cup V^{1,a'}$ is of the form (\ref{v2}), (\ref{va}) or (\ref{vb}).

\medskip
{\bf Case} $|V^{1,a}|=|V^{1,a'}|=2$ and $d=6$.  By Lemma \ref{rr3}$ii$ and \ref{ci} we may assume that for every $v\in V^{1,a}\cup V^{1,a'}$ every cover $C_v$ of the word $v$ by words from $W$  contains a code $P^{1,b}\cup P^{1,b'}$, where $q$-equivalent codes $P^{1,b}_B$ and $P^{1,b'}_B$ for $B=\{4,5,6\}$ and $q=bbb$, are given in Table 2 at the positions 1-16  (that is, $P^{1,b}_{B}=K, P^{1,b'}_{B}=M$, where $K,M$ are given in that table). Recall that, by Lemma \ref{cody}$iv$, $w_2=w_3=b$ for every $w\in P^{1,b}\cup P^{1,b'}$. (Compare also Subsection 2.6.)

We use algorithm {\small{\sc CoverCode}} in which we take $U=V^{1,a}\cup V^{1,a'}$ and $P_i=P^{1,b}\cup P^{1,b'}$ for every $i\in [4]$ (recall that $V^{1,a}\cup V^{1,a'}$ is of the form (\ref{v2})). For $v\in V^{1,a}\cup V^{1,a'}$ let $\ka C_v$ denote the set of all twin pair free covers $C_v\subset S^d$ of $v$ such that $|C_v|\leq 11$, $P^{1,b}\cup P^{1,b'}\subset C_v$ and $|C^{1,l}_v|\leq 3$ for $l\in \{a,a'\}$ (the last inequality steams from  the condition $|W^{1,a}|,|W^{1,a'}|\leq 3$). We have $|\ka C_v|=9233$ for $v\in \{aabbbb,a'babbb\}$ and $|\ka C_v|=2214$ for $v\in \{aa'abbb,a'aa'bbb\}$. (There is no cover of the words $aa'abbb$ and $a'aa'bbb$ with eleven elements and containing $P^{1,b}_{B}=K$ and $P^{1,b'}_{B}=M$, where $K,M$ are given in the last row of Table 2.)  

\smallskip
In the result of the computations we obtained seven non-empty covers $C_1,...,C_7$ of the code $V^{1,a}\cup V^{1,a'}$ such that each of them has less than seventeen words and satisfies the condition  $|W^{1,a}|, |W^{1,a'}|\leq 3$. 

Every code $C_i$, $i\in [6]$, has thirteen words and $D_k(C_i)=((0,1),(10,2),(0,0))$ for some $k\in [6]$ and every $i\in [6]$. In the similar way as in the proof of the second part of ($i$) in Lemma \ref{rr3} we show that   if $W$ is a twin pair free disjoint equivalent code to $V$ and $C_i\subset W$ for $i\in [6]$, then $D_k(W)=((2,2),(10,2),(0,0))$. Thus, $|W|=16$. 

The computations show that for every $i\in [6]$ and for every three words $v,u,w\in S^6$, if $W=C_i\cup \{v,u,w\}$ is a twin pair code which satisfies the condition $|W^{1,a}|, |W^{1,a'}|\leq 3$ and the conclusions of Lemma \ref{rr5} applied to $W^{i,l},W^{i,l'}$ for $i\in [d]$ and $l\in \{a,b\}$, then there are $j\in [6]$ and $l\in S$ such that $|W^{j,l}|=5$ and $|W^{j,l'}|=1$ or $|W^{j,l}|=6$ and $|W^{j,l'}|=1$ and in both cases $W^{j,l'}_{j^c}\not\sqsubseteq W^{j,l}_{j^c}$. This, by Lemma \ref{po1m}, means that $W^{j,l}\sqsubseteq V^{j,l}$ and $W^{j,l'}\sqsubseteq V^{j,l'}$. Then, by Lemma \ref{po7}, $|V^{j,l}|\geq 10$ and $|V^{j,l'}|\geq 5$, and consequently, by Lemma \ref{rr2}, $V=V^{j,l}\cup V^{j,l'}$. By Lemma \ref{r24}, $|V|\geq 24$, a contradiction.

The last code $C_7$ has the form:
$$
C_7^{1,a}\cup C_7^{1,a'}=\{ab'bbbb,\;aa'b'bbb,\;a'b'a'bbb,\;a'bb'bbb\}
$$
and 
$$
C_7^{1,b}\cup C_7^{1,b'}=\{bbbaaa,\;bbba'a'a',\;bbbbaa',\;bbba'ba,\;bbbaa'b,\;b'bbbbb\}.
$$

Now we show that every twin pair free code $W$ such that $C_7\subset W$ and $|W| \in \{13,...,16\}$ cannot be equivalent to $V$. Suppose on the contrary that there is a twin pair free code $W$ with $C_7\subset W$ which is equivalent to $V$.

Let $w^1=aa'b'bbb$, $v^1=aa'abbb$ and $w^2=a'bb'bbb$, $v^2=a'babbb$. We have $v^1,v^2 \in  V^{1,a}\cup V^{1,a'}$ and $w^1,w^2 \in C_7\subset W$. Since $w^1_{3^c}=v^1_{3^c}$ and $w^1_3=b',v^1_3=a$, by ({\bf P}), $v^1_{3^c}\sqsubseteq V^{3,a'}_{3^c}$. For the same reason, $v^2_{3^c}\sqsubseteq V^{3,a'}_{3^c}$. Thus, $v^1_{3^c}, v^2_{3^c}\sqsubseteq V^{3,a'}_{3^c}$. Therefore, by Lemma \ref{po7}$ii$, $|V^{3,a'}|\geq 7$. Moreover, $D_3(C_7)=((0,1),(7,2),(0,0))$. Note that, by Theorem \ref{t12} and (${\bf P}$), $V^{3,b}\neq\emptyset$. Thus, by Lemma \ref{rr2}, the sets $V^{3,l}, W^{3,l}$ are nonempty for $l\in \{a,a',b,b'\}$ and, by Lemma \ref{cy}, $|V^{3,a'}|+|V^{3,b'}|+|V^{3,c'}|=|W^{3,a'}|+|W^{3,b'}|+|W^{3,c'}|$ and $|V^{3,a}|+|V^{3,b}|+|V^{3,c}|=|W^{3,a}|+|W^{3,b}|+|W^{3,c}|$. Since  $|V^{3,a'}|+|V^{3,b'}|\geq 8$  and  $|W^{3,a}|+|W^{3,b}|\geq 8$, we obtain  $V^{3,l}, W^{3,l}=\emptyset$  for $l\in \{c,c'\}$, $|W^{3,b}|=|V^{3,a'}|=7$ and $|W^{3,a}|=|V^{3,b'}|=1$. Hence, $|W|=16$. We shall show that
\begin{equation}
\label{ro1}
D_3(W)=((1,1),(7,7),(0,0))\;\; {\rm and}\;\;  D_3(V)=((7,7),(1,1),(0,0)).
\end{equation}
Suppose on the contrary that $(|W^{3,a}|,|W^{3,a'}|)\neq (1,1)$. Since $|W^{3,a}|=1$, by Lemma \ref{rr3}$i$, we have $|W^{3,a'}|\geq 5$. Clearly, there is $v\in V\setminus (V^{3,a}\cup V^{3,a'})$ such that $\breve{v}\cap \bigcup E(W^{3,a}\cup W^{3,a'})\neq\emptyset$, otherwise $W^{3,a'}\sqsubseteq V^{3,a'}$ and then, by Lemma \ref{po7}$iii$, $|V^{3,a'}|>7$, a contradiction. By Lemma \ref{po1m}, $W^{3,a}_{3^c}\sqsubseteq W^{3,a'}_{3^c}$, and then, by Lemma \ref{rr1}, $V^{3,a}_{3^c}\sqsubseteq V^{3,a'}_{3^c}$. Therefore, by Lemma \ref{po7}, $|V^{3,a}|\leq 2$. Consequently, it is easy to check that $|W^{3,a'}|=6$ and $|V^{3,a}|=2$. Thus,  $|W^{3,a}|=1$, $|W^{3,a'}|=6$. We shall show that this is not possible. Since $\breve{v}\cap \bigcup E(W^{3,a}\cup W^{3,a'})\neq\emptyset$, the codes $U^{3,a}_{3^c}$, $U^{3,a'}_{3^c}$ are, by Lemma \ref{wac}, $v_{3^c}$-equivalent, where $U^{3,l}=\{w\in W^{3,l}\colon \breve{v}\cap \breve{w}\neq\emptyset\}$ for $l\in \{a,a'\}$. If $|U^{3,a'}|=5$, then, by Lemma \ref{cody}$i$, \ref{cody}$iv$  and \ref{ci} the set $P^{3,l}=\{u\in V\setminus (V^{3,a}\cup V^{3,a'})\colon \breve{u}\cap \bigcup E(U^{3,a}\cup U^{3,a'})\neq\emptyset\}$ contains two words for every $l\in \{b,b'\}$. A contradiction, because $P^{3,b'}\subset V^{3,b'}$ and $|V^{3,b'}|=1$.  If $|U^{3,a'}|=6$, then, again by Lemma \ref{cody}$i$, \ref{cody}$iv$ and \ref{ci}, there is a twin pair in $V$ or in $W$, which is not true. This completes the proof of the equalities (\ref{ro1}). 

\smallskip
Clearly, we have
$$
W^{3,b}=C_7^{3,b},\;\;\;\; C_7^{3,b'}\subset W^{3,b'}.
$$
Moreover, $m(W^{3,l}_{3^c}\setminus W^{3,l'}_{3^c})\leq 1$ for $l\in \{b,b'\}$, by (\ref{mb}) and (\ref{ro1}). 

Now we show that $|W^{1,c}\cup W^{1,c'}|\geq 4$. If it is not true, by Lemma \ref{rr3}$i$, $|W^{1,c}|=|W^{1,c'}|=1$, and then $D_1(W)=((2,2),(5,5),(1,1))$ which gives, by Lemma \ref{wa} and \ref{po7}$i$, 
$D_1(V)=((2,2),(1,1),(5,5))$ (compare the beginning of the proof of Lemma \ref{ogr}).
Note now that $w_{1^c}\sqsubseteq W^{1,b}_{1^c}$, where $w=b'bbbbb$. Since $w\in W^{1,b'}$, by (\ref{5n5}), there is no twin pair free code $W^{1,b}\cup W^{1,b'}$ such that $|W^{1,l}|=5$ for $l\in \{b,b'\}$ and  the equality (\ref{dif}) is valid, where $|V^{1,l}|=1$ for $l\in \{b,b'\}$. A contradiction, by Lemma \ref{rr1}. 

Since $|W^{1,c}\cup W^{1,c'}|\geq 4$, it follows that at least two words from the set $W^{1,c}\cup W^{1,c'}$ have to belong to the set $W^{3,b'}$. But it is easy to compute that for every two words $w,u\in S^d\setminus C_7^{3,b'}$ such that $w_1,u_1\in \{c,c'\}$ if $Q=C_7^{3,b'}\cup \{w,u\}$, where recall $C_7^{3,b'}\subset W^{3,b'}$, is a twin pair free code, then $m(W^{3,b}_{3^c}\setminus Q_{3^c})>4$. 
This shows that the inequality $m(W^{3,b}_{3^c}\setminus W^{3,b'}_{3^c})\leq 1$ does not hold. Thus, (\ref{ro1}) is not possible, and consequently $V$ cannot contain a code $V^{1,l}\cup V^{1,l'}$ such that $|V^{1,l}|=|V^{1,l'}|=2$ and $|W^{1,l}|,|W^{1,l'}|\leq 3$ for $l\in S$ and $d=6$.

\medskip
{\bf Case} $|V^{1,a}|=|V^{1,a'}|=2$ and $d=5$. The code $V^{1,a}\cup V^{1,a'}$ is obtained from (\ref{v2}) by skipping the last letter $b$ in every word of  (\ref{v2}). By Lemma \ref{rr3}$ii$, \ref{rr3}$iv$  and \ref{ci} we may assume that for every $v\in V^{1,a}\cup V^{1,a'}$ every cover $C_v$ of the word $v$ by words from $W$  contains the code 
$P=\{bbbab,bbba'a,b'bbba,b'bbaa'\}.$ Moreover, by Lemma \ref{ogr},  $|W^{1,l}|\leq 3$ for $l\in \{a,a'\}$ and $|C_v|\leq 11$ for every $v\in V^{1,a}\cup V^{1,a'}$. We use algorithm {\small{\sc CoverCode}} in which we take $U=V^{1,a}\cup V^{1,a'}$ and $P_i=P$ for every $i\in [4]$. Our aim is to compute all covers $C_U\subseteq W$ of $U$ with the above properties. 

We have $|\ka C_v|=2159$ for $v\in \{aabbb,a'babb\}$ and $|\ka C_v|=1130$ for $v\in \{aa'abb,a'aa'bb\}$. 

The computations show that every such cover $C_U$ of $U$ which satisfies the condition $|C_U^{1,l}|\leq 3$ for $l\in \{a,a'\}$
has more than sixteen words. Thus, $V$ cannot contain a code $V^{1,l}\cup V^{1,l'}$ such that $|V^{1,l}|=|V^{1,l'}|=2$ and $|W^{1,l}|,|W^{1,l'}|\leq 3$ for $l\in S$ and $d=5$.

\medskip
{\bf Case} $|V^{1,a}|=2,\;|V^{1,a'}|=3$ and $d=6$. Recall that the code $V^{1,a}\cup V^{1,a'}$ has the form (\ref{va}) or (\ref{vb}).
For the same reasons as in the previous case we compute all twin pair free covers of $V^{1,a}\cup V^{1,a'}$ such that for $v\in V^{1,a}\cup V^{1,a'}$ every cover $C_v$ of the word $v$ by words from $W$ contains the  code $P=\{bbbbab, bbbba'a, b'bbbba, b'bbbaa'\}$, $|C_v|\leq 11$ and $|W^{1,a}|, |W^{1,a'}|\leq 3$. We use algorithm {\small{\sc CoverCode}} in which we take $U=V^{1,a}\cup V^{1,a'}$ and $P_i=P$ for every $i\in [5]$.

We have $|\ka C_v|=6799$ for $v\in \{aabbbb\}$, $|\ka C_v|=1946$ for $v\in \{a'ba'a'bb, a'ababb\}$ and $|\ka C_v|=222$ for $v\in \{aa'a'a'bb,a'aaa'bb\}$ if $V^{1,a}\cup V^{1,a'}$ is of the form (\ref{va}). If it is of the form (\ref{vb}), then $|\ka C_v|=1946$ for $v\in \{aaabbb,\; aba'a'bb,a'aba'bb\}$ and $|\ka C_v|=222$ for $v\in \{a'aaabb,\; a'a'a'a'bb\}$.

\smallskip
The computations show that every such cover $C_U$ of $U$ which satisfies the condition $|C_U^{1,l}|\leq 3$ for $l\in \{a,a'\}$ has more than sixteen words. Thus, $V$ cannot contain a code $V^{1,l}\cup V^{1,l'}$ such that $|V^{1,l}|=2, |V^{1,l'}|=3$ and $|W^{1,l}|,|W^{1,l'}|\leq 3$ for $l\in S$.

\medskip
{\bf Case} $|V^{1,a}|=2,\;|V^{1,a'}|=3$ and $d=5$. Skipping the last letter in (\ref{va}) and (\ref{vb}) we obtain the set $V^{1,a}\cup V^{1,a'}$ for $d=5$. Since, by Lemma  \ref{cody}$iv$, $w_2=w_3=w_4=b$ for every $w\in W\setminus (W^{1,a}\cup W^{1,a'})$ such that $\breve{w}\cap \bigcup E(V^{1,a}\cup V^{1,a'})\neq\emptyset$, by Lemma \ref{ci}, there is a twin pair in $V$ or in $W$. A contradiction.

This, together with the conclusions of the previous cases, contradicts Lemma \ref{ogr}. The proof of the lemma is completed.
\hfill{$\square$}

\section{The proof of Theorem \ref{keli1}}

A code $U\subset S^d$ with $2^d$ words is called a {\it partition code}. Any realization $f(U)$ of a partition code $U$ is a minimal partition (compare Subsection 2.3).

The main theorem of the paper steams from the following

\begin{tw}
\label{ke2}
There are no two disjoint equivalent polybox codes $V,W\subset S^d$ without twin pairs for $d\leq 6$, where $|V|\in \{13,...,16\}$. In particular, if $U\subset S^d$, $d\leq 7$, is a partition code such that there is $i\in [d]$ for which the set $U^{i,l}$ is non-empty for every $l\in \{a,b,c,d\}\subset S$, then $U$ contains  a twin pair.
\end{tw}
\proof
To prove the first part of the theorem assume on the contrary that $V$ and $W$ are disjoint twin pair free equivalent polybox codes, where $|V|\in \{13,...,16\}$. 

Let $d=4$. We give the proof for the case $|V|=13$. It was computed in \cite{CS2} that every code $Q\subseteq \{a,a',b,b'\}^4$ having more than twelve words contains a twin pair.  It follows from Theorem 2 of \cite{DIP}, that there are three words $v,w$ and $u$ such that $U=V\cup \{v,w,u\}$ is a partition code. Clearly, there are $i\in [4]$ and $l\in S$ such that $U^{i,l}=V^{i,l}$ and $U^{i,s}\neq\emptyset$ for $s\in \{a,a',b,b',c,c'\}\subseteq S$. Since $|U^{i,l}|\leq 6$, the code $U^{i,l}$ is rigid (see the claim in the proof of Lemma \ref{rr2}). Thus, if $U^{i,l'}\setminus \{v,w,u\}\neq\emptyset$, then $(U^{i,l'}\setminus \{v,w,u\})_{i^c}\cap U^{i,l}_{i^c}\neq\emptyset$, and consequently there is a twin pair in $V$. A contradiction. If   $U^{i,l'}\setminus \{v,w,u\}=\emptyset$, then  $U^{i,s}=V^{i,s}$ and  $U^{i,s'}=V^{i,s'}$ for some $s\in S\setminus \{l,l'\}$. Since $U^{i,s}$ is rigid (we have  $|U^{i,s}|\leq 6$), there is a twin pair in  $U^{i,s}\cup U^{i,s'}$, a contradiction.

In the case $d\in \{5,6\}$ we shall made computations for $d=6$ and $d=5$ separately. By Lemma \ref{3li} we may assume that $V,W\subset S^d$, where $S=\{a,a',b,b'\}$.

Let $d=6$. By Corollary \ref{le12}, we assume that every word $v\in V$ is covered by at most ten words from $W$.  We may assume that the word $v=bbbbbb$ belongs to $V$ (compare Subsection 2.6) and let $\ka C^{n}_v$ be the family of all twin pair free covers of the word $v$ having $n$ words for $n\in \{5,...,10\}$. The family $\ka C_6(v)=\bigcup _{n=5}^{10}\ka C^n_v$ contains 2058920 covers, and 104 of them are non-isomorphic codes. We denote the set of such codes by $\ka N_6$. 


\smallskip
Now we shall proceed as follows: For every cover $U\in \ka N_6$ of $v$, based on algorithm {\small {\sc CoverCode}}, we compute all twin pair free covers $C_U$ of the code $U$ which satisfy the conclusions of Lemma \ref{rr5}, $C_U\cap U=\emptyset$, $v\in C_U$ and $|C_U|\leq 16$. 
In the next step we compute all twin pair free covers $C_{C_U}$ of $C_U$ which satisfy the conclusions of Lemma \ref{rr5}, $C_{C_U}\cap C_U=\emptyset$, $U\subset C_{C_U}$ and $|C_{C_U}|\leq 16$. 

For $U\in \ka N_6$ by $\ka C_U$ we denote the family of all covers $C_U$, and $\ka C^2_U$ stands for the family of all $C_{C_U}$. As as we shall see $\ka C^2_U=\emptyset$ for every $U\in \ka N_6$. This implies 
that $V$ and $W$ cannot be equivalent. Indeed, if $V$ and $W$ are equivalent then we may assume that there is $U\in \ka N_6$ such that $U\subset W$ (compare Subsection 2.6) and there is a cover $C_U\subset V$ of $U$ such that the cover $C_{C_U}=\{w\in W\colon \breve{w}\cap \bigcup E(C_U)\neq\emptyset\}$ has at most sixteen words (because $|W|\leq 16$), $U\subset C_{C_U}$ (because $U\subset W$), $C_{C_U}\cap C_U=\emptyset$ (because $V\cap W=\emptyset$) and $C_{C_U}$ satisfies the conclusions of Lemma \ref{rr5} (because $W$ satisfies them). Thus, $\ka C^2_U\neq \emptyset$.


\smallskip
To compute $C_U$, $U\in \ka N_6$, using algorithm {\small {\sc CoverCode}} we need also initial configurations $P_i$, $i\in [|U|]$ described in this algorithm. Let $U=\{u^1,...,u^n\}$. In our case $P_i=\{v\}$ for every $u^i\in U$, where recall $v=bbbbbb$. This means that every word $u^i\in U$ is covered by codes $C_{u^i}\in \ka C_{u^i,P_i}$, where $\ka C_{u^i,P_i}$ is the family of all codes $C_{u^i}\in \ka C_6(u^i)$ such that $P_i\subset C_{u^i}$.

\smallskip
To indicate codes $P_i$ for the computations of a cover $C_{C_U}$ of the code $C_U$, $C_U\in \ka C_U$, let $C_U=\{v^1,...,v^k\}$ and $P_i=\{u\in U\colon \breve{u}\cap \breve{v}^i\neq\emptyset\}$ for $i\in [k]$.  Clearly, if $C_{v^i}\subset C_{C_U}$ is a cover of $v^i$, then $P_i\subset C_{v^i}$. Thus, computing  $C_{C_U}$ for every $i\in [k]$ we use only covers from the family $\ka C_{v^i,P_i}$ consisting of all codes $C_{v^i}\subset \ka C_6(v^i)$ such that $P_i\subset C_{v^i}$.

\medskip
As the result of the computations we obtained $\ka C_U=\emptyset$ for all, but five, codes $U\in \ka N_6$. These five families $\ka C_U$ contains 42, 48, 48, 24 and 24 covers $C_U$, and the corresponding codes $U$ contain 5, 6, 7, 8 and 9 words, respectively. Computing in all five cases the sets $\ka C^2_U$ we obtain each time  $\ka C^2_U=\emptyset$. This completes the proof of the first part of the theorem for $d=6$.

\smallskip
In the same way we prove the theorem for $d=5$. In this case, by Corollary \ref{le12}, it is assumed that every word $V$ is covered by at most eleven words from $W$. As for $d=6$ we may assume that the word $v=bbbbb$ belongs to $V$
The family $\ka C_5(v)=\bigcup _{n=5}^{11}\ka C^n_v$ contains 738680 covers, and 232 of them are non-isomorphic codes. The set of these non-isomorphic codes is denoted by $\ka N_5$. 

\smallskip
The computations show that $\ka C_U=\emptyset$ for all, but two, codes $U\in \ka N_5$. These two families $\ka C_U$ contains 324 and 8 covers $C_{U}$, and the corresponding codes $U$ contain five and seven words, respectively. Computing in these two cases the sets $\ka C^2_U$ we obtain each time  $\ka C^2_U=\emptyset$. This completes the proof of the first part of the theorem for $d=5$.

\medskip
To prove the second part of the theorem, assume that there is $l\in \{a,b,c,d\}$ such that $|U^{i,l}|\leq 12$. Then, as we showed in Theorem 29 of  \cite{Kis22}, there is a twin pair in $U$. Therefore, we may assume that $|U^{i,l}|\geq 13$ for every $l\in \{a,b,c,d\}$. Then there is $l\in \{a,b,c,d\}$ such that $|U^{i,l}|\in \{13,...,16\}$. The codes $U^{i,l}_{i^c}$ and $U^{i,l'}_{i^c}$ are equivalent, and hence, by the first part of the theorem, $U^{i,l}_{i^c}\cap U^{i,l'}_{i^c}\neq\emptyset$ or there is a twin pair in $U^{i,l}_{i^c}$ or in $U^{i,l'}_{i^c}$. In both case the set $U^{i,l}\cup U^{i,l'}$ contains a twin pair.
\hfill{$\square$}

\medskip
We are ready to prove the main theorem of the paper.

\medskip
{\it Proof of Theorem \ref{keli1}}. Since $r^+(T)=4$, it follows that there are $x\in \er^7$ and $i\in [7]$ such that there are sets $A^1,...,A^4\subset [x_i,1+x_i]$ with $A^j\not \in \{A^k,(A^k)^c\}$ for $k,j\in [4], k\neq j$, where $(A^k)^c=[x_i,1+x_i]\setminus A^k$ and $\ka F_x=\ka F_x^{i,A^1}\cup \ka F_x^{i,(A^1)^c}\cup \cdots \cup \ka F_x^{i,A^4}\cup \ka F_x^{i,(A^4)^c}$ (see Subsection 2.2).

Thus, we may write down a partition code $U$, whose realization is the minimal partition $\ka F_x$ in the alphabet $S=\{a,a',b,b',c,c',d,d'\}$, where $U^{i,l}\neq\emptyset$ for above mentioned $i\in [7]$ and every $l\in S$. (A way of receiving of such code from a minimal partition is given at the end of  Subsection 2.7 of \cite{Kis22}). From Theorem \ref{ke2} we infer that $U$ contains a twin pair. Then $\ka F_x$ contains such a pair, and consequently there is a twin pair in $[0,1)^7+T$.
\hfill{$\square$} 

\medskip
Theorem 1.1 together with Theorem 1.1 of \cite{Kis} and \cite{Kis22} show that Keller's conjecture is true for cube tilings $[0,1)^7+T$ with $r^+(T)\geq 4$. On the other hand the result of Debroni et al. (\cite{De}) shows that it is true for cube tilings $[0,1)^7+T$ with $r^-(T)\leq 2$ (\cite[Introduction]{Kis}). Thus, we obtain

\begin{wn}
\label{kon3} 
If $[0,1)^7+T$ is a counterexample to Keller's conjecture, then $r^+(T)=3$.
\hfill{$\square$}
\end{wn}

\medskip
Theorem \ref{ke2} provides also the following

\smallskip
{\it Proof of Theorem \ref{clic}}. It is straightforward: If $V$ is a clique with $|V|=2^d$, then $V$ is a partition codes without twin pairs. This, by Theorem \ref{ke2}, is impossible. Thus, $|V|<2^d$.
\hfill{$\square$}

\medskip
{\it Proof of Corollary \ref{wk}}. If $[0,1)^7+T$ is a counterexample, then by Corollary \ref{kon3}, there is $x \in \er^7$ such that a partition code $U$ of the minimal partition $\ka F_x$ can be  written down in the alphabet $S=\{a,a',b,b',c,c'\}$. Clearly, $V=U$ is a clique in the Keller graph on $\{a,a',b,b',c,c'\}^7$ with $|V|=128$, as it is a partition code. Since $r^+(T)=3$, there is $i\in [7]$ such that $V^{i,a}\neq\emptyset$, $V^{i,b}\neq\emptyset$ and $V^{i,c}\neq\emptyset$. If $0<|V^{j,l}|\leq 16$ for some $j\in [7]$ and some $l\in S$, then the similar arguments as those used at the end of the proof of Theorem \ref{ke2} show that $V$ contains a twin pair. Then the tiling $[0,1)^7+T$ contains a twin pair, a contradiction.  

If $V$ is a clique with 128 elements, then for $i\in [7]$ let $f_i(a)=[0,1)+2\zet,f_i(a')=[1,2)+2\zet,$ $f_i(b)=[\frac{1}{3},\frac{4}{3})+2\zet,f_i(b')=[\frac{4}{3},\frac{7}{3})+2\zet$ and $f_i(c)=[\frac{1}{2},\frac{3}{2})+2\zet,f_i(b')=[\frac{3}{2},\frac{5}{2})+2\zet$. The realization $f(V)$ of $V$ (see Subsection 2.3) is a twin pair free minimal partition of the $7$-box $\er^7$ (see Subsection 2.1). It is easy to see that at the same time it can be viewed as a $2$-periodic cube tiling of $\er^7$ without tiwn pair.      
\hfill{$\square$}

\medskip
Finally, we can give a new proof of Keller's conjecture in dimensions $d\leq 6$ (recall that it was originally proved by Perron in 1940; another proof is given in \cite{Kis}).  

\begin{tw}
\label{per}
In every cube tiling $[0,1)^d+T$ of $\er^d$, where $d\leq 6$, there is a twin pair.
\end{tw}
\proof
It is enough to show that in every partition code $U\subseteq S^d$, $d\leq 6$, there is a twin pair. It is obvious for a partition code $U$ such that $U=U^{i,l_i}\cup U^{i,l'_i}$ for every $i\in [d]$, where $l_i\in S$ for $i\in [d]$ (as it is isomorphic to the binary code $\{0,1\}^d$). Let $U\neq  U^{i,l}\cup U^{i,l'}$ for some $i\in [d]$ and every $l\in S$. Then there is $s\in S$ such that $0<|U^{i,s}|\leq 16$. This, as we have seen in the proof of Theorem \ref{ke2}, implies that $U$ contains a twin pair.
\hfill{$\square$}





\begin{thebibliography}{AB}
\bibitem{CS1} K. Corr\'adi and S. Szab\'o,  Cube tiling and covering a complete graph, {\it Discrete Math.} 85 (1990), 319--321.
\bibitem{CS2} K. Corr\'adi and S. Szab\'o,  A combinatorial approach for Keller's conjecture, {\it Period. Math. Hungar.} 21 (1990), 95--100.
\bibitem{De} J. Debroni, J.D Eblen, M.A. Langston, W. Myrvold, P. Shor and D. Weerapurage, A complete resolution of the Keller maximum clique problem, {\it Proceedings of the Twenty-Second Annual ACM-SIAM Symposium on Discrete Algorithms}, 2011.
\bibitem {DIP} M. Dutour Sikiric, Y. Itoh and A. Poyarkov, Cube packings, second moment and holes, \textit{European J. Combin.} {\bf 28} (2007), 715--725.
\bibitem {GKP} J. Grytczuk, A. P. Kisielewicz and K. Przes{\l}awski, Minimal Partitions of a Box into Boxes, {\it Combinatorica} 24 (2004), 605--614.
\bibitem {Ke1} O. H. Keller, \"Uber die l\"uckenlose Erf\"ullung des Raumes
mit W\"urfeln,  {\it J. Reine Angew. Math.} 163 (1930), 231--248. 
\bibitem {Kis} A. P. Kisielewicz, Rigid polyboxes and Keller's conjecture, to appear in {\it Adv. of Geom.}, available on arXiv:1304.1639
\bibitem {Kis22} A. P. Kisielewicz, M. {\L}ysakowska, On Keller's conjecture in Dimension Seven, {\it Electron. J. Combin.}  22 (2015), \#P1.16, pp. 44. 
\bibitem {KP} A. P. Kisielewicz, K. Przes{\l}awski, Polyboxes, cube tilings and rigidity, {\it Discrete Comput. Geom.} 40 (2008), 1--30.
\bibitem {LS1} J. C. Lagarias and P. W. Shor, Keller's cube-tiling
conjecture is false in  high dimensions, {\it Bull. Amer. Math. Soc.} 27 (1992), 279--287.
\bibitem {LS2} J. C. Lagarias and P. W. Shor, Cube tilings of $\er^d$ and nonlinear codes,
{\it Discrete Comput. Geom.}  11 (1994), 359--391.
\bibitem {La} J. Lawrence, Tiling $\er^d$ by translates of orthants, Convexity and Related Comb. Geometry Proc. of the Second Univ. of Oklahoma Conf. (1982), 203--207.
\bibitem {LP1} M. {\L}ysakowska and K. Przes{\l}awski, Keller's conjecture on the existence of columns in cube tilings of $\er^n$,   {\it Adv. Geom.} 12  (2012), 329--352.
\bibitem {M} {J. Mackey}, A cube tiling of dimension eight with no facesharing, {\it Discrete
Comput. Geom.} 28 (2002), 275--279.
\bibitem {P} {O. Perron}, \"Uber l\"uckenlose Ausf\"ullung des $n$-dimensionalen Raumes durch
kongruente W\"urfeln, {\it Math. Z.} 46 (1940), 1--26.
\bibitem {Sz2} S. Szab\'o, A reduction of Keller's conjecture, {\it Period. Math. Hungar.} 17 (1986),  265--277.

\end{thebibliography}
\end{document}